\documentclass[12pt]{article}

\usepackage[letterpaper,margin =1in]{geometry}
\usepackage{amssymb}
\usepackage{amsmath}
\usepackage{amsthm}
\usepackage{mathtools}
\usepackage{enumerate}
\usepackage{fancyhdr}
\usepackage{braket}
\usepackage{titlesec}
\usepackage{xcolor}
\usepackage{hyperref}
\usepackage{lmodern} 

\usepackage{tikz-cd} 
\usetikzlibrary{babel}

\usepackage{tikz}
\usetikzlibrary{shadows, arrows}

\usepackage{scalerel,stackengine}
\stackMath
\newcommand\reallywidehat[1]{%
\savestack{\tmpbox}{\stretchto{%
  \scaleto{%
    \scalerel*[\widthof{\ensuremath{#1}}]{\kern-.6pt\bigwedge\kern-.6pt}%
    {\rule[-\textheight/2]{1ex}{\textheight}}
  }{\textheight}%
}{0.5ex}}%
\stackon[1pt]{#1}{\tmpbox}%
}

\usepackage{todonotes}

\definecolor{skyblue}{rgb}{0.85,0.85,1}

\titleformat{\section}
  {\normalfont\Large\bfseries}{\thesection}{1em}{}

\titleformat{\subsection}
  {\normalfont\large\bfseries}{\thesubsection}{1em}{}

\DeclareFontFamily{OT1}{pzc}{}
\DeclareFontShape{OT1}{pzc}{m}{it}%
             {<-> s * [1.200] pzcmi7t}{}
      
\newlength\tindent
\DeclareMathAlphabet{\mathscr}{OT1}{pzc}{m}{it}

\providecommand{\del}{\partial}
\providecommand{\eps}{\varepsilon}

\newenvironment{frcseries}{\fontfamily{frc}\selectfont}{}

\newcommand{\ve}{\varepsilon}

\newcommand{\sgn}{\textrm{sgn}}
\newcommand{\re}{\textrm{Re}}

\newcommand{\cA}{\mathcal{A}}

\newcommand{\C}{\mathbb{C}}

\newtheorem{thr}{Theorem}[section]

\newtheorem{coro}[thr]{Corollary}
\newtheorem{prop}[thr]{Proposition}

\theoremstyle{definition}
\newtheorem{defi1}[thr]{Definition}

\theoremstyle{remk}
\newtheorem{remark1}[thr]{Remark}
\newenvironment{remk}{\begin{remark1}\rm}{\hfill $\triangle$ \end{remark1}}

\makeatletter
\newcommand{\tpitchfork}{%
  \vbox{
    \baselineskip\z@skip
    \lineskip-.52ex
    \lineskiplimit\maxdimen
    \m@th
    \ialign{##\crcr\hidewidth\smash{$-$}\hidewidth\crcr$\pitchfork$\crcr}
  }%
}
\makeatother

\begin{document}

\title{Nonlinear stability of shock-fronted travelling waves under nonlocal regularization}
\author{Ian Lizarraga\thanks{School of Mathematics and Statistics, The University of Sydney, {\tt ian.lizarraga@sydney.edu.au} (corresponding author)}  \: and Robert Marangell\thanks{School of Mathematics and Statistics, The University of Sydney, {\tt robert.marangell@sydney.edu.au}.}}
\maketitle

\begin{abstract}
We determine the nonlinear stability of shock-fronted travelling waves arising in a reaction-nonlinear diffusion PDE, subject to a fourth-order spatial derivative term multiplied by a small parameter $\varepsilon$ that models {\it nonlocal regularization}. Motivated by the authors' recent stability analysis of shock-fronted travelling waves under viscous relaxation, our numerical analysis is guided by the observation that there is a fast-slow decomposition of the associated eigenvalue problem for the linearised operator. In particular, we observe an astonishing reduction of the complex four-dimensional eigenvalue problem into a {\it real} one-dimensional problem defined along the slow manifolds; i.e. slow eigenvalues defined near the tails of the shock-fronted wave for $\varepsilon = 0$ govern the point spectrum of the linearised operator when $0 < \varepsilon \ll 1$. 
\end{abstract}


\section{Introduction}

{\it Reaction-nonlinear diffusion partial differential equations} (hereafter RND PDEs) have recently been attracting interest as a natural mathematical model for {\it coherent propagation with sharp fronts}. On the one hand, they can be derived as continuum limits of discrete models of biological competition and exclusion \cite{simpson1,simpson2}; on the other, they are now known to admit shock-fronted travelling wave solutions \cite{li,li2}.\\

  At the same time, the technique of smoothing shock-type solutions in nonlinear diffusion models via high-order regularization is by now well-motivated and studied; a prototypical example in this context is the Cahn-Hilliard model (see e.g. Pego and Penrose's work in \cite{pego1989}). Regularized nonlinear diffusion models with advection have also been considered by Witelski in \cite{witelski1995,witelski1996}. An interesting feature of high-order regularizations is that they give rise to distinct rules for shock selection. In \cite{witelski1996}, Witelski used singular perturbation theory to derive generalised area rules for shock selection in nonlinear diffusion models, corresponding to an interpolated family of viscous and nonlocal regularizations. This result has recently been extended to the RND PDE setting; see \cite{proceedings22}.\\

We consider the following RND PDE with a fourth-order spatial derivative term modelling nonlocal regularization:\footnote{We follow the convention from our earlier work \cite{lizarraga}: barred functions of $(x,t)$ are used for the travelling wave dynamics in phase space, and unbarred functions of $(x,t)$ are used for linearized quantities.}
\begin{align} \label{eq:master}
\frac{\del \bar{U}}{\del t} &= \frac{\del}{\del x} \left( D(\bar{U}) \frac{\del \bar{U}}{\del x} \right) + R(\bar{U}) - \eps^2 \frac{\del^4 \bar{U}}{\del x^4}\end{align}

for $(x,t) \in \mathbb{R}\times \mathbb{R}$. The small parameter $\eps \geq 0$ characterizes the strength of the regularizing effect. We suppose that the potential function $F(\bar{U}) := \int D(\bar{U}) d\bar{U}$  is nonmonotone, i.e. $D(\bar{U}) < 0$ within an interval $(a,b) \subset (0,1)$, and we also suppose that the reaction term is bistable, being pinned at $U = 0$ and $U = 1$. \\

Using geometric singular perturbation theory (GSPT), Li and his coauthors \cite{li} showed the existence of a one-parameter family of (smoothed) shock-fronted travelling waves for system \eqref{eq:master}. Using a coordinate frame that follows the travelling wave, the PDE \eqref{eq:master} can be written as a closed four-dimensional system of ODEs in two fast and two slow variables. The shock-fronted travelling waves are then found as heteroclinic orbits connecting a saddle-type fixed point at $\bar{U} = 1$ to one at $\bar{U} = 0$. These fixed points lie on disjoint branches of two-dimensional saddle-type slow manifolds, and the shock segment corresponds to a fast jump between these manifolds. Li et al. also derive the well-known {\it equal area rule} for these waves by utilizing the Hamiltonian structure of the (fast) layer subsystem.\\

The objective of this paper is to determine the {\it nonlinear stability} of this family of shock-fronted travelling waves. Following the standard strategy, we calculate the spectrum $\sigma(\mathcal{L})$ of a corresponding sectorial linearized operator $\mathcal{L}$ associated with the PDE \eqref{eq:master}. Our task is to show that except for a simple translational eigenvalue that necessarily exists at the origin, the entire spectrum is bounded within the left-half complex plane in a `nice' way, as we now describe. The total spectrum is decomposed into its continuous and point components as $\sigma(\mathcal{L}) = \sigma_e(\mathcal{L}) \cup \sigma_p(\mathcal{L})$. Our first result concerns the location of the continuous spectrum and its implications for nonlinear stability: completing the analysis that we initiated in the Concluding Remarks of \cite{lizarraga}, we show in Sec. \ref{sec:eigenvalueproblem} that the linearised operator is {\it sectorial}.  We can already highlight a major difference with the viscous relaxation limit; there, the continuous spectrum was similarly shown to be bounded inside the left-half complex plane, but it was asymptotically vertical, obstructing sectoriality of the operator. \\

It then remains to show that the translational eigenvalue is simple and that the rest of the point spectrum is bounded within the left-half complex plane, and this is the primary focus of this paper. Following standard geometric stability theory for travelling waves (see e.g. \cite{kapprom} and \cite{bjornstability}), we characterize values $\lambda \in \sigma_p (\mathcal{L})$ in the point spectrum via bifurcations of the eigenvalue problem $(\mathcal{L}-\lambda )v = 0$  posed on the underlying linearized subspaces along the travelling wave.  Specifically, we are interested in tracking the complex two-plane subbundle of unstable subspaces along travelling wave solutions of \eqref{eq:master}, beginning at the asymptotic steady-state $\bar{U} = 1$ of \eqref{eq:master}. An eigenvalue $\lambda \in \sigma_p(\mathcal{L})$ is found whenever this bundle has a nontrivial intersection with the corresponding stable subbundle continued from the stable subspace at $\bar{U} = 0$.\\

Numerical methods to locate these intersections focus on robustly approximating a suitable {\it Evans function}, whose roots coincide with values in the point spectrum; we refer the reader to some standard references \cite{kapprom,bjornstability}. A recent technique that has found much success in numerical stability calculations is the computation of the winding number of a  {\it Riccati-Evans function} along a suitably chosen contour in the complex plane; see \cite{harley,ledoux2} and references therein. Unlike a traditional Evans function, the Riccati-Evans function is typically not analytic but at best meromorphic, and is defined by considering the induced flow of the linear eigenvalue problem on a suitable (chart of the) Grassmannian (in the present case, we use $\text{Gr}(2,4)$, the Grassmannian of complex 2-planes in $\mathbb{C}^4$). As demonstrated by Harley et al. in \cite{harley}, blowups of the induced matrix Riccati equation (i.e. coordinate pole singularities) can be moved by selecting charts judiciously in the domain of interest.\\

Our main result (see Sec. \ref{sec:riccatievans}) is the use of a Riccati-Evans function to calculate the point spectrum. We provide strong evidence that the point spectrum consists of only a simple translational eigenvalue at the origin, and another simple real eigenvalue of negative real part. We conclude that the family of shock-fronted travelling waves is nonlinearly stable.\\

An equally significant objective of this paper is to explain the output of the Riccati-Evans calculation by demonstrating an unambiguous {\it fast-slow splitting} of the eigenvalue problem. We couch our result within the geometric framework developed by Alexander, Gardner, and Jones in the series \cite{AGJ,GJ,jones}.  In recent work by the authors \cite{lizarraga}, we showed rigorously that for $\eps > 0$ sufficiently small, the eigenvalue problem in the (three-dimensional) viscous relaxation case was controlled by a slow eigenvalue problem defined on the critical manifolds for $\eps = 0$, together with a newly-defined {\it jump map}. Note that the unstable subbundle in the viscous relaxation case is a line bundle (after a suitable time reversal).  In the present case, we demonstrate numerically that the unstable 2-plane bundle splits in such a way that a `fast' line subbundle always makes a trivial unstable-to-unstable connection for $\eps > 0$ sufficiently small, hence contributing no eigenvalues. On the other hand, a `slow' line subbundle continues to be controlled by a slow eigenvalue problem. \\
 
One of the advantages of producing this splitting is that the analysis of the slow problem in the nonlocal case carries over {\it identically} from that of the viscous case: indeed, the reduced two-dimensional linearized system defined on the critical manifold is exactly the same. For example, we can state immediately that there are no slow eigenvalues with nonzero imaginary part (see Sec. 8 in \cite{lizarraga}). By studying the linear flow on projective space, we therefore obtain the promised reduction of the original eigenvalue problem to a real one-dimensional problem.  \\

We wish to contrast our results with the travelling wave stability analysis of another four-dimensional singularly-perturbed problem with two fast and two slow variables considered by Gardner and Jones in \cite{GJ}. In particular, we report two crucial differences. The first is that the `fast unstable-to-unstable' connections that we find do not persist in the singular limit, instead requiring a very delicate break in the (linearized) dynamics near a heteroclinic orbit in the intermediate fast layer as soon as $\eps > 0$ (see Sec. \ref{sec:fastconnections}). The second difference is the requirement of a nontrivial jump map to define the slow problem, which can be clearly observed in our results (see Sec. \ref{sec:sloweigs}). These two aspects are related, as shown in the earlier setting considered by Gardner and Jones: specifically, an elephant trunk lemma enforces a trivial jump map for the slow dynamics (see Theorem 5.3(b) in \cite{GJ}). However, our fast eigenvalue problem {\it degenerates} in the singular limit $\eps \to 0$. As a consequence, we cannot construct a suitable `fast elephant trunk' that Gardner and Jones use in their analysis to track the fast unstable solution of the eigenvalue problem through the shock layer. \\

 The geometric issue of splitting the unstable 2-plane subbundle is also intimately related to the analytic problem of factorizing an Evans function into fast and slow reduced components  over the entire travelling wave (see \cite{derijk}). Revisiting the present problem using the analytic factorization procedure proposed by de Rijk et al. is likely to be fruitful. \\

The rest of this paper is organized as follows. In Sec. \ref{sec:existence}, we review the construction of the one-parameter family of shock-fronted travelling waves in the PDE \eqref{eq:master}. In Sec. \ref{sec:eigenvalueproblem} we write down the spatial eigenvalue problems for the linearized operator of \eqref{eq:master} and remind the reader of previous results concerning the essential spectrum. In Sec. \ref{sec:riccatievans} we define the Riccati-Evans function and compute its roots, thereby locating the point spectrum for the singular perturbation parameter $\eps = 10^{-4}$. In Sec. \ref{sec:fastconnections} we describe numerical observations of the fast (un)stable-to-(un)stable line bundle connections in the eigenvalue problem, and in Sec. \ref{sec:sloweigs} we demonstrate that the slow line subbundle is well-approximated by a reduced eigenvalue problem defined near the slow manifolds of the system, together with a jump map. In this section we also describe how the slow eigenvalue problem defined in the singular limit $\eps \to 0$ controls the generation of eigenvalues in the point spectrum in the `full' problem when $\eps > 0$. We conclude with some open questions in Sec. \ref{sec:conclusion}.

\section{Existence of shock-fronted travelling waves} \label{sec:existence}

We begin by reviewing how Li and his coauthors \cite{li} used GSPT to construct a family of shock-fronted travelling waves of \eqref{eq:master} near the singular limit $\eps = 0$ (i.e. `small' regularisation). We freely use well-known definitions from {\it Fenichel theory} (see \cite{fenichel} and the comprehensive introductions in \cite{joneslec},\cite{kuehn}) when referring to concepts like normal hyperbolicity and locally invariant slow manifolds.

We define the diffusion and reaction terms in \eqref{eq:master} concretely as
\begin{align} 
D(\bar{U}) &= 6(\bar{U}-7/12)(\bar{U}-5/6) \label{eq:diffusionterm}\\
R(\bar{U})  &= 5 \bar{U} (1-\bar{U})(\bar{U}-1/5). \label{eq:reaction}
\end{align}

These definitions are selected for consistency with \cite{li} and \cite{lizarraga}. We use the standard technique of expressing the PDE \eqref{eq:master} in terms of the frame $(\zeta,t') = (x-ct,t)$, where $c$ is an unknown constant parameterizing the wavespeed, and then seeking steady-states of the resulting transformed PDE. By using a Li\'enard representation for the resulting closed four-dimensional system of ODEs (see Sec. 2.1 in \cite{lizarraga}), we obtain the following  system in two fast variables $(\bar{U},\bar{W})$ and two slow variables $(\bar{P},\bar{V})$:
\begin{equation} \label{eq:slowphase}
\begin{aligned}
\eps\dot{\bar{U}} &= \bar{W}\\
\eps\dot{\bar{W}} &=  F(\bar{U})-\bar{V}\\
\dot{\bar{P}} &= -R(\bar{U})\\
\dot{\bar{V}} &= \bar{P}-c\bar{U}.
\end{aligned}
\end{equation}

Here $F(\bar{U}) :=\int D(\bar{U})$ denotes an integral of $D(\bar{U})$, which we fix as
\begin{align} \label{eq:potentialdef}
F(\bar{U}) &= 2\bar{U}^3 - 4\bar{U}^2 + \frac{21}{8}\bar{U}.
\end{align}

In terms of the stretched travelling wave variable $\xi := \zeta/\eps$, we obtain the following system which is equivalent to \eqref{eq:slowphase} when $\eps > 0$:\footnote{In keeping notational consistency with \cite{GJ,lizarraga} we use the lowercase alphabet for quantities defined w.r.t. the stretched variable $\xi$ and uppercase otherwise, i.e. $a(\xi) := A(\zeta) = A(\eps \xi)$ for an arbitrary function $A(\zeta)$.}
\begin{equation} \label{eq:fastphase}
\begin{aligned}
\bar{u}' &= \bar{w}\\
\bar{w}' &=  F(\bar{u})-\bar{v}\\
\bar{p}' &= -\eps R(\bar{u})\\
\bar{v}'&= \eps(\bar{p}-c\bar{u}).
\end{aligned}
\end{equation}

Note that the dot notation denotes the derivative with respect to $\zeta$ and the prime notation denotes the derivative with respect to $\xi$. The system \eqref{eq:slowphase} (or \eqref{eq:fastphase}) admits three fixed points corresponding to zeroes of the bistable reaction term $R(\bar{u})$. The fixed points $q_- = (0,0,0,0)$ and $q_+ = (1,0,c,F(1))$ are saddle-type for each $\eps > 0$ and $c > 0$, both of which have two-dimensional local stable and unstable manifolds. The travelling waves that we seek are now realized as heteroclinic orbits connecting the saddle point at $\bar{u} = 1$ (via a weak unstable direction) to the one at $\bar{u} = 0$ (via a weak stable direction). \\

The singular limit $\eps \to 0 $ of systems \eqref{eq:slowphase} and \eqref{eq:fastphase} are distinct. The system \eqref{eq:fastphase} limits to the fast {\it layer problem}, defined as

\begin{equation} \label{eq:layerproblem}
\begin{aligned}
\bar{u}' &= \bar{w}\\
\bar{w}' &=  F(\bar{u})-\bar{v}\\
\bar{p}' &=  0\\
\bar{v}' &=  0.
\end{aligned}
\end{equation}

\begin{figure}[t!] 
\centering
\includegraphics[width=1.05\textwidth]{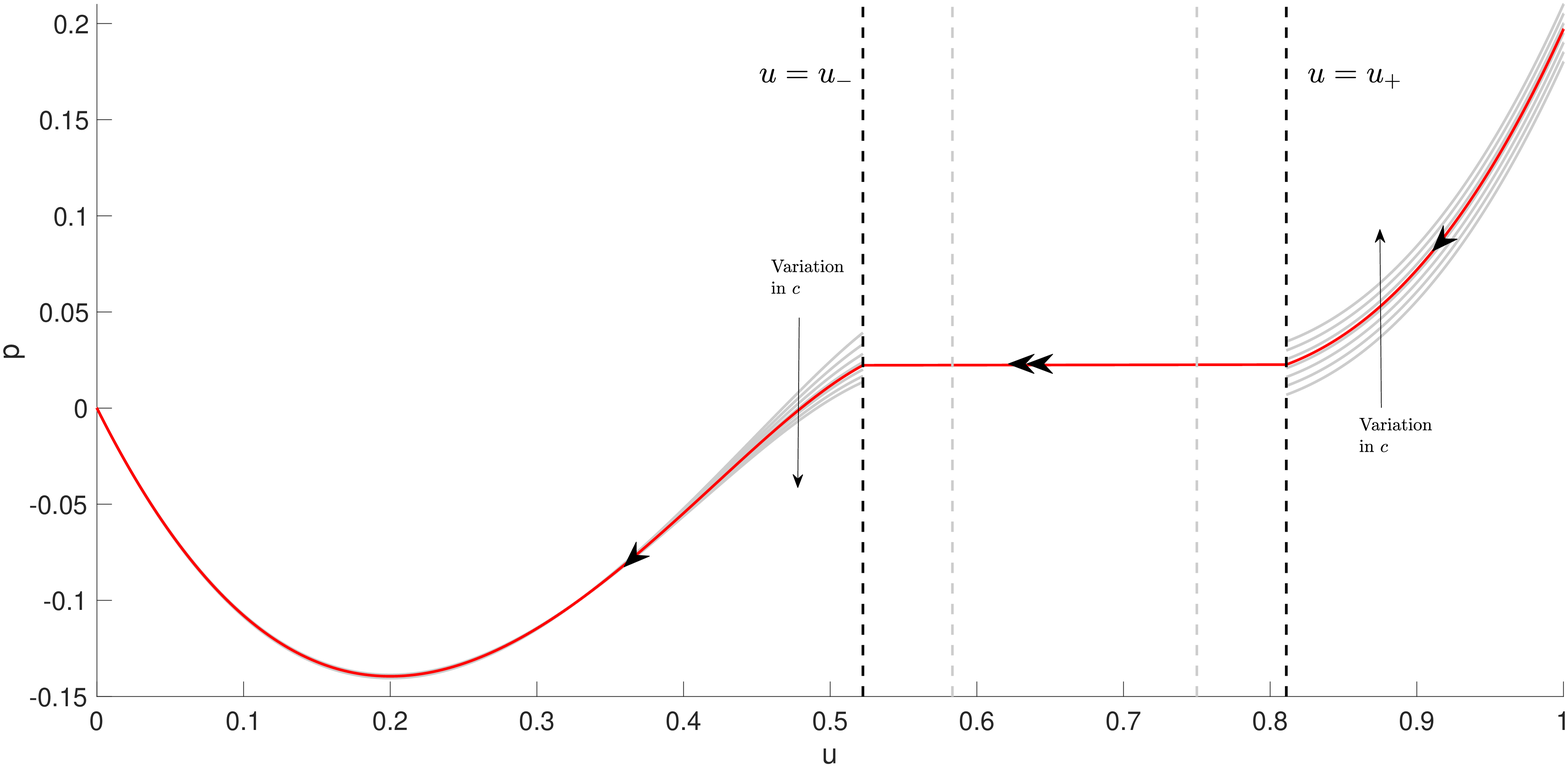}
\caption{Segments of $W^s(\bar{Q}_-)$ and $W^u(\bar{Q}_+)$ arising from the slow flow \eqref{eq:reduced2}, as $c$ is varied within the range $[0.17,0.22]$. The singular homoclinic orbit $\Gamma_0$ is depicted by the red curve for $c_0 \approx 0.1968109995$, and the corresponding fast shock is also overlaid. The fold lines at $\bar{u} = \bar{u}_{f,\pm}$ are depicted as gray dashed lines, and the $\bar{u}$-values $\bar{u} = \bar{u}_{\pm}$ satisfying the equal area rule are depicted as black dashed lines.}
\label{fig:transversehet}
\end{figure}

The slow variables $(\bar{p},\bar{v})$ can be interpreted as parameters in the singular limit. On the stretched `timescale,' we obtain a two-dimensional problem containing a two-dimensional cubic manifold of equilibria called the {\it critical manifold}:
\begin{align*}
S &:= \{(\bar{u},\bar{w},\bar{p},\bar{v}) \in \mathbb{R}^4: \bar{v} = F(\bar{u}),~\bar{w} = 0\}.
\end{align*}
The layer problem \eqref{eq:layerproblem} defines the dynamics along the {\it fast fibers} away from this critical set. We can classify the {\it normal hyperbolicity} along points of $S$ by studying the eigenvalues of the linearization of \eqref{eq:layerproblem}. There are always two zero eigenvalues, called the {\it trivial eigenvalues}, corresponding to the two-dimensional tangent spaces at every point in $S$. The signs of the remaining two {\it nontrivial} eigenvalues  at a point $x \in S$ determines the normal hyperbolicity at that point. A calculation of the Jacobian shows us that $S$ is a {\it saddle-type} critical manifold on the sets $S_- := S \cap \{\bar{u} < u_{f,-} := 3/4\}$ and $S_+ := S \cap \{\bar{u} > \bar{u}_{f,+} := 5/6\}$,  and is {\it repelling} when $\bar{u}_{f,-} < \bar{u} < \bar{u}_{f,+}$. The critical manifold contains two {\it fold lines} at $\bar{u} = \bar{u}_{f,-}$ and $\bar{u} = \bar{u}_{f,+}$ along which it loses normal hyperbolicity via a simple zero eigenvalue crossing. Note that the fold lines occur at roots of the diffusion $D(\bar{u})$. The travelling waves that we consider in this paper stay far away from the fold lines, in contrast to the viscous case \cite{lizarraga}. \\

As highlighted in \cite{li}, the layer problem \eqref{eq:layerproblem} has the structure of a Hamiltonian system with
\begin{align*}
H(u,w) &= \frac{1}{2}w^2 - G(u) + vu,
\end{align*}

where $G(u)$ is any integral of $F(u)$. Fast connections adjoining the outer saddle-type critical branches exist when there are heteroclinic orbits of \eqref{eq:layerproblem} connecting saddle points at values $u = u_- < u_{f,-}$ to those at $u = u_+ > u_{f,+}$. Using the fact that these are level sets of the Hamiltonian, we obtain the {\it equal area rule} 
\begin{align*}
\int_{u_-}^{u_+} (v-F(u))\,du &= 0
\end{align*}
which constrains such pairs of values $u_{\pm}$ . In view of the symmetry of the graph of the cubic potential function \eqref{eq:potentialdef} about its inflection point at $\bar{u} = 2/3$, we find the exact expressions $$\bar{u}_{\pm} = \frac{1}{12}(8\pm \sqrt{3}).$$

\begin{figure}[t!] 
\centering
(a)\includegraphics[width=0.85\textwidth]{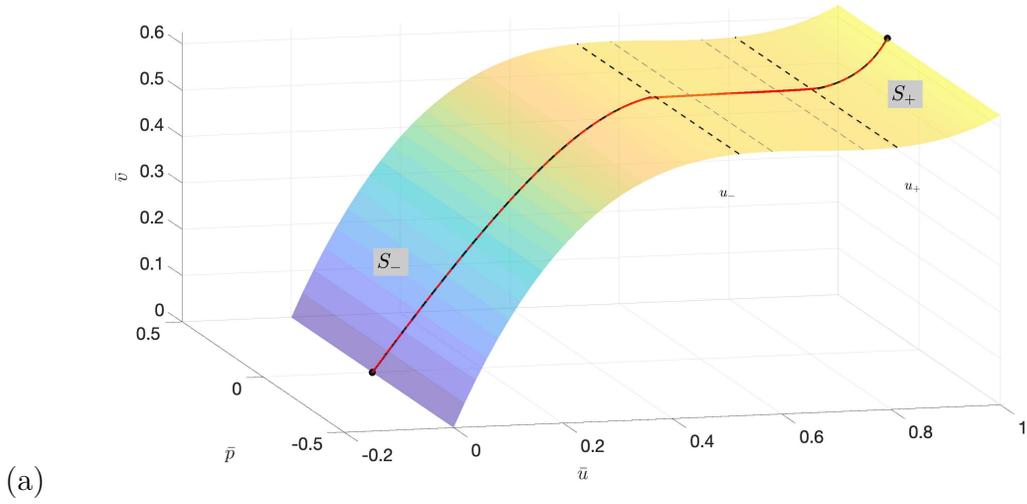}\\
(b)\includegraphics[width=0.85\textwidth]{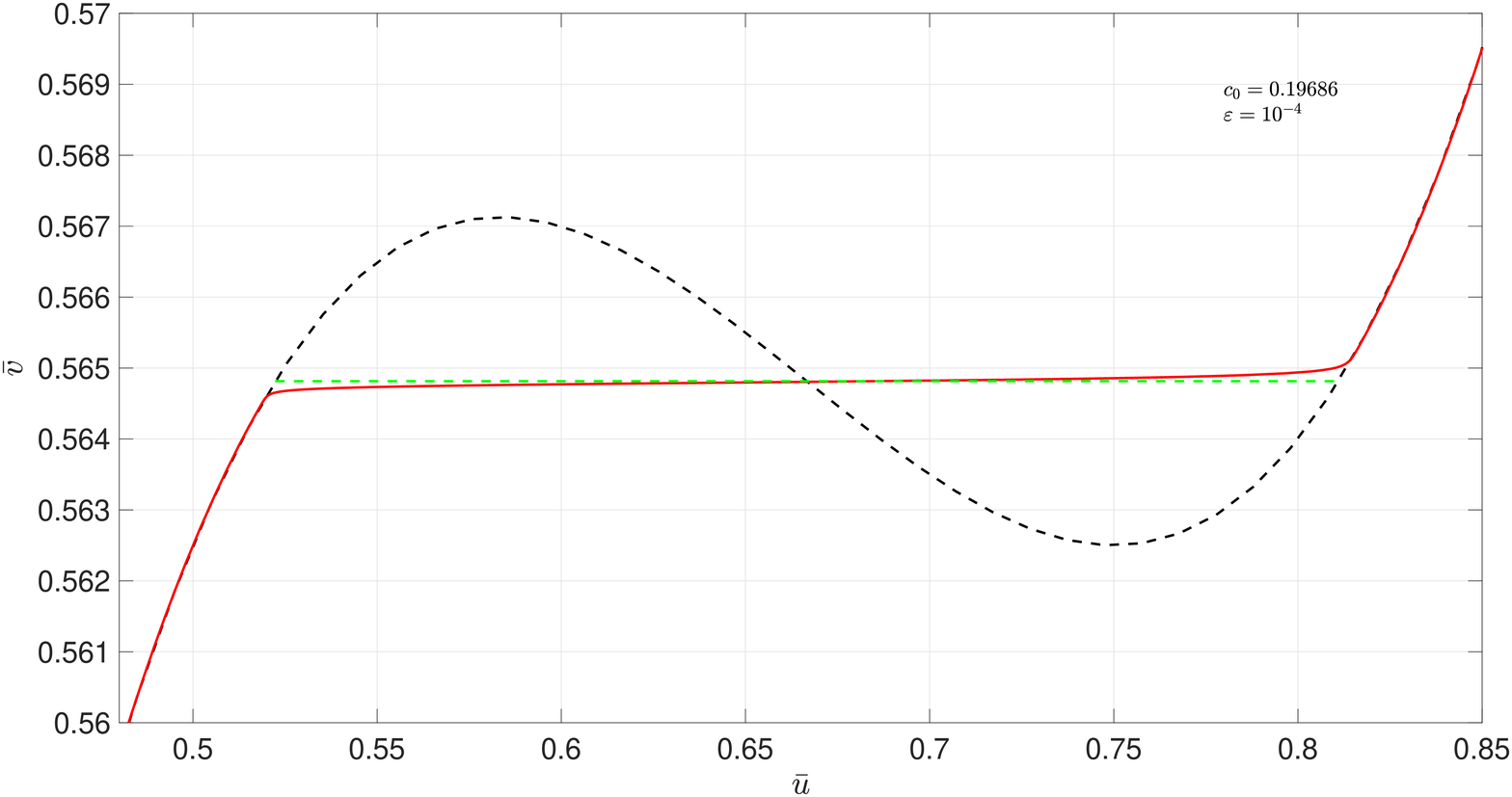}
\caption{(a) Heteroclinic orbit of \eqref{eq:slowphase} for the parameter set $(\eps,c) = (10^{-4},0.19686)$ on a $(\bar{u},\bar{p},\bar{v})$-projection of the four-dimensional phase space (red curve) together with segments of the singular homoclinic orbit on the critical manifold $S$ (black dashed curves). The saddle-type branches $S_{\pm}$ are also identified. (b) Magnified view near the folds of $S$. The green dashed line depicts the equal area rule and the black dashed line depicts the critical manifold. The orbit was computed using the \texttt{bvp4c} boundary-value solver in MATLAB R2021a with a relative error tolerance of $10^{-8}$. }
\label{fig:fullhet}
\end{figure}

We now write down the singular limit of \eqref{eq:slowphase}. As $\eps \to 0$, we obtain the two-dimensional {\it reduced problem}, constrained to the critical manifold $S$:
\begin{equation} \label{eq:reducedproblem}
\begin{aligned}
\bar{V} &= F(\bar{U})\\
\bar{W} &=0 \\
\dot{\bar{P}} &= -R(\bar{U})\\
\dot{\bar{V}} &= \bar{P}-c\bar{U}.
\end{aligned}
\end{equation}

The reduced problem defines the slow dynamics on the critical manifold in the singular limit. The constraint $\dot{\bar{V}} = D(\bar{U}) \dot{\bar{U}}$, which specifies that the reduced flow should be everywhere tangent to $S$, provides a convenient closed representation of the reduced problem as
\begin{equation} \label{eq:reduced2}
\begin{aligned}
\dot{\bar{U}} &= \frac{1}{D(\bar{U})}(\bar{P} - c\bar{U})\\
\dot{\bar{P}} &= -R(\bar{U}).
\end{aligned}
\end{equation}

Note that the saddle singularities $\bar{q}_-$ and $\bar{q}_+$ for $\eps >0$ limit to saddle-type fixed points of \eqref{eq:reduced2} at $\bar{Q}_- = (0,0)$ resp. $\bar{Q}_+ = (1,c)$.  Solution segments of the fast and slow reduced problems \eqref{eq:layerproblem} and \eqref{eq:reduced2} can then be formally concatenated with the prospect of finding {\it singular heteroclinic orbits} connecting $Q_-$ and $Q_+$.\\

 In Fig. \ref{fig:transversehet}, we report the existence of such an orbit $\Gamma_0$, which is also transverse with respect to variation in the wavespeed parameter $c$. One of the important achievements of GSPT is the use of Fenichel theory together with other estimates (such as the exchange lemma and blow-up theory, if necessary) to rigorously construct heteroclinic orbits as robust perturbations of these singular objects. 
We refer the reader to \cite{li} for the detailed construction of a one-parameter family of heteroclinic orbits $\{\Gamma_{\eps}\}_{\eps \in (0,\bar{\eps}]}$ for the full system \eqref{eq:slowphase}. Here $\bar{\eps}$ denotes a sufficiently small upper bound on $\eps$. A member of this family is plotted in Fig. \ref{fig:fullhet}.

\section{Point and essential spectra} \label{sec:eigenvalueproblem}

We now turn to the stability analysis of the family of heteroclinic orbits $\{\Gamma_{\eps}\}_{\eps \in (0,\bar{\eps}]}$. For given $\eps>0$, we represent the corresponding member of this family by the solution \\$\bar{X}(\zeta,\eps) = (\bar{U}(\zeta,\eps),\bar{W}(\zeta,\eps),\bar{P}(\zeta,\eps),\bar{V}(\zeta,\eps))$ and consider perturbed solutions of \eqref{eq:master} of the form $\tilde{U}(\zeta,t) = \bar{U}(\zeta,\eps)+\delta e^{\lambda t}U(\zeta)+\mathcal{O}(\delta^2)$. At linear order we find the temporal eigenvalue problem in the linearised operator $\mathcal{L}$, defined as:
\begin{align} \label{eq:linearisedoperator}
\lambda U &= cU_{\zeta} + (D(\bar{U}) U)_{\zeta\zeta} + R'(\bar{U})U_{\zeta} - \eps^2 U_{\zeta \zeta \zeta \zeta} =: \mathcal{L}(U)
\end{align}

 After collecting terms of linear order in $\delta$ and choosing a Li\'{e}nard coordinate representation for $\mathcal{L}$ that is compatible with that of the travelling wave phase space dynamics, we arrive at the nonautonomous linear dynamical system
\begin{equation} \label{eq:linslow}
\begin{aligned}
\eps\dot{U} &= W\\
\eps\dot{W} &=  D(\bar{U})U - V\\
\dot{P}&= (\lambda-R'(\bar{U}))U\\
\dot{V} &= P - cU
\end{aligned}
\end{equation}

on the (slow) $\zeta$ timescale and

\begin{equation} \label{eq:linfast}
\begin{aligned}
u' &= w\\
w' &=  D(\bar{u})u - v\\
p' &= \eps (\lambda-R'(\bar{u}))u\\
v' &= \eps(p - cu)
\end{aligned}
\end{equation}

on the stretched (fast) $\xi$ timescale.\\

 Let us denote the linear system \eqref{eq:linslow} as $Z' = A(\xi,\lambda,\eps)Z$. We define the {\it asymptotically constant matrices} associated with $A(\xi,\lambda,\eps)$, given by $$A_{\pm}(\lambda,\eps):=\lim_{\xi\to \pm \infty}A(\xi,\lambda,\eps)$$ (note: we analogously define the matrices $a_{\pm}$ in the natural way). \\
 
Now we want to compute the spectrum of the linearised operator $\mathcal{L}$. We begin with the essential spectrum, $\sigma_e(\mathcal{L})$. Via Weyl's essential spectrum theorem (see, for example \cite{kapprom}), we can characterize the essential spectrum of $\mathcal{L}$ by determining the values of $\lambda$ for which the signatures of $A_{\pm}(\lambda,\eps)$ differ; equivalently, the {\it Fredholm borders} which bound the essential spectrum are characterized by the condition that one of $A_-(\lambda,\eps)$ and $A_+(\lambda,\eps)$ has an imaginary (spatial) eigenvalue $ik$.   The calculations are given in the  Concluding Remarks of \cite{lizarraga} (the nonlocal case we consider is the case $a = 0$ in that computation). Summarising, we have that the dispersion relations forming the Fredholm borders are given by: 

\begin{equation}\label{eq:dispersion}
\lambda_\pm = - \ve^2 k^4 - D(\bar{U}^\pm)k^2 + R'(\bar{U}^\pm) + i c k 
\end{equation}
where 
\begin{equation}
\bar{U}^{\pm} := \lim_{\zeta \to \pm \infty} \bar{U} 
\end{equation}
 (and analogously for the lower case letters). \\
 
These are seen to partition $\C$ into five disjoint sectors, distinguished by the signatures of $A_\pm(\lambda, \ve)$. We call region to the right of the essential spectrum, where both $A_\pm(\lambda,\ve)$ have signature $(2,2)$, $\Omega$. That is 

$$
\Omega := \left\{ \lambda \in \C \bigg{|} \sgn(A_-) = \sgn(A_+) = (-,-,+,+) \right\}.
$$

\begin{figure}[t!] 
\centering
\includegraphics[width=0.85\textwidth]{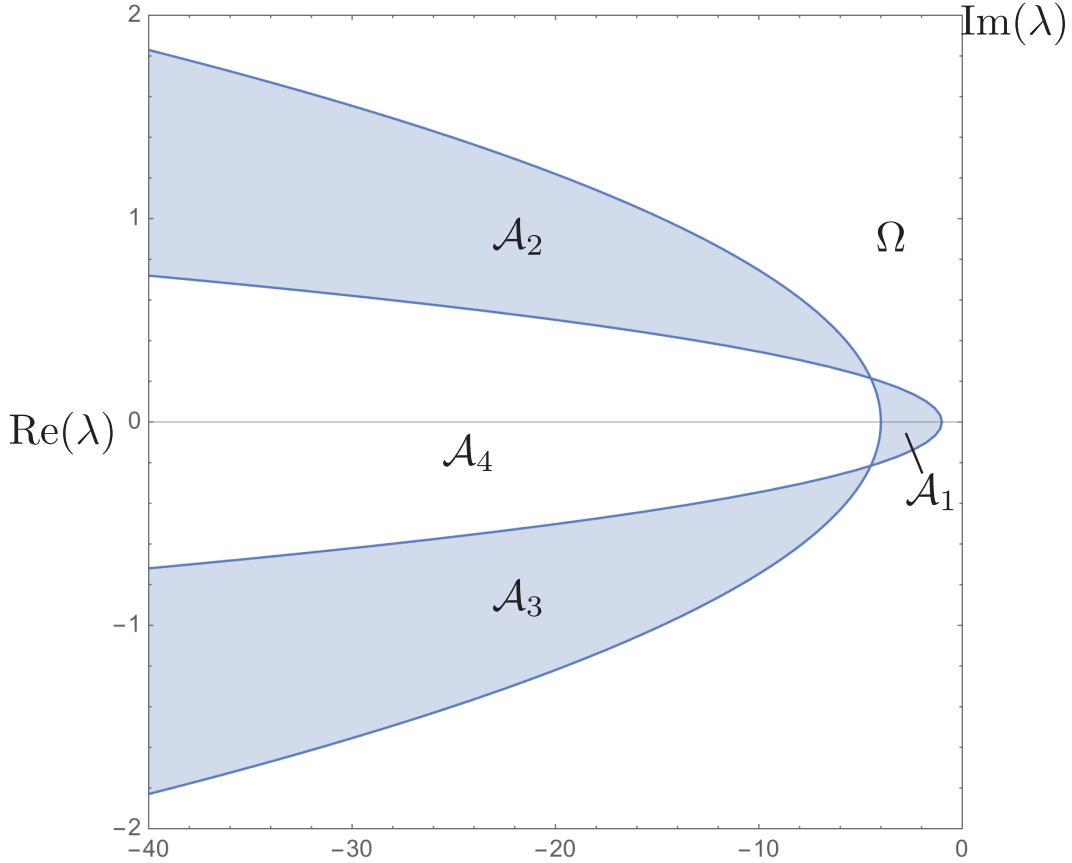}\\
\caption{A plot of the essential spectrum (shaded region) of the operator $\mathcal{L}$ in \eqref{eq:linearisedoperator}. The five regions, $\Omega$ and $\mathcal{A}_i$ with $i = 1,2,3,4$ are shown. The bounding lines are the Fredholm borders given by the dispersion relations in \eqref{eq:dispersion}. For the figure, $\ve = 10^{-5}$, and $c =.19681$, with $D$ and $R$ as in \eqref{eq:diffusionterm} and \eqref{eq:reaction}.
}
\label{fig:essential-pic4}
\end{figure}

We use $\mathcal{A}_i$ with $i = 1,2,3,4 $ to denote the remaining regions. We note that the essential spectrum is contained within the half plane $\{x+iy \in \mathbb{C}: x\leq -1= R'(\bar{U}^-)\}$, and the behavior for large $|\lambda|$ is dominated by the quartic term, $ -\ve^2k^4$ (see Figure \ref{fig:essential-pic4} and Table \ref{tab:4dsigs}).\\


\begin{table}
\begin{centering}
\begin{tabular}{c|| c | c }
Region & $\sgn(A_{-}(\lambda,\ve))$ & $\sgn(A_{+}(\lambda,\ve)f)$ \\ 
\hline
$\Omega$ & $(-, -, +,+)$ & $(-,-,+,+)$ \\ 
$\cA_{1}$ & $(-, -,-, +)$ & $(-,-,+,+)$ \\
$\cA_{2}$ & $(-, -, + +)$ & $(-,-,-,+)$ \\
$\cA_{3}$ & $(-, -, +, +)$ & $(-,-,-,+)$ \\
$\cA_{4}$ & $(-, -,-, +)$ & $(-,-,-,+)$ \\
\end{tabular}
\caption{The signatures of the asymptotically constant matrices which define the regions $\Omega$ and $\mathcal{A}_i$.}
\label{tab:4dsigs}
\end{centering}
\end{table}

 For $\lambda \in \Omega$, the asymptotic operators $A_{\pm}(\lambda,\eps)$ are hyperbolic, having eigenvalues $\mu^{\pm}_i$ (with $i=1,2,3,4$) that fall within the hierarchy $$\text{Re}(\mu^{\pm}_4)  \ll \text{Re}(\mu^{\pm}_3) < 0 <  \text{Re}(\mu^{\pm}_2) \ll \text{Re}(\mu^{\pm}_1).$$ This estimate can be shown to hold uniformly along the segments of the wave near to the slow manifolds $S_{\eps}$, i.e. by studying the linearisations of $A(\xi,\lambda,\eps)$ defined on the left and right saddle branches of the slow manifold for $|\xi|$ sufficiently large, the corresponding eigenvalues $\mu^{\pm}_j(\xi)$ inherit this hierarchy. See e.g. \cite{GJ,lizarraga} for examples of such estimates in three- and four-dimensional problems.\\
 
   Note that this hierarchy implies the existence of {\it fast unstable} and {\it fast stable} directions for the eigenvalue problem corresponding to the eigenspaces of $\mu^{\pm}_1$ and $\mu^{\pm}_4$, respectively, together with a $2$-plane {\it slow subbundle}  given by the span of eigenspaces corresponding to $\mu^{\pm}_{2,3}$. \\

We seek to determine whether the family of traveling waves $\Gamma(c(\eps),\eps)$ is {\it nonlinearly stable}, i.e. whether for fixed $\eps > 0$, every solution $U(t)$ of \eqref{eq:master} with initial data sufficiently near to the wave will tend to a translate of the wave at an exponential rate as $t \to \infty$. General stability theory (see e.g. \cite{henry}, as well as Prop. 2.2 and Cor. 2.1 in \cite{AGJ}) provides broad hypotheses for which it is sufficient to check the {\it linearised} stability of the wave. In fact, the operator $\mathcal{L}$ is sectorial, as we shall see, so it is sufficient to check the spectral stability of the wave. We then have linearised and nonlinear stability of the wave due to standard results in the theory (see for example \cite{AGJ,henry,kapprom}). \\

The dispersion relations in \eqref{eq:dispersion} show that the essential spectrum is contained in the left half plane. We now give a geometric characterization of the point spectrum. We restrict our attention to the set of eigenvalue parameters $\lambda \in \Omega$, the subset of the complex plane to the right of the Fredholm borders of the essential spectrum.  For each such $\lambda \in \Omega$, it is a direct calculation to check that for $\eps >0$ sufficiently small, the asymptotically constant matrices $A_{\pm}(\lambda,\eps)$ each have two eigenvalues of positive real part and two eigenvalues of negative real part (in particular, note that this is true for $\lambda = 0$, as indicated in the discussion of the hyperbolic fixed points in \ref{sec:existence}). By general theory (see \cite{AGJ}), we may extend an {\it unstable 2-plane bundle} $\varphi^{+}(\zeta,\lambda,\eps)$ from the saddle point point $q_+$ at $\bar{u} = 1$, using the linearised flow \eqref{eq:linslow}: it is the unique bundle that converges to $\text{span}\{v_1^+,v_2^+\}$ as $\zeta \to -\infty$, where $v^+_{1,2}$ span the eigenspaces corresponding to the unstable eigenvalues  $\mu^{+}_{1,2}$ of \eqref{eq:linslow} at the steady-state $\bar{U} = 1$. There is also an analogously defined {\it stable 2-plane bundle} $\varphi^{-}(\zeta,\lambda,\eps)$ extending outward from the stable eigendirections of the fixed point $q_-$ perched at $\bar{u} = 0$. \\

The general theory (see \cite{AGJ}) tells us that $\lambda \in \sigma_p (\mathcal{L})$ if and only if $\varphi^-$ and $\varphi^+$ have an intersection at some common value $\zeta = \zeta_0$ (and hence for each $\zeta \in \mathbb{R}$) . Defining the section 
\begin{align} \label{eq:sigmasection}
\Sigma &= \{\bar{u} = 2/3\}
\end{align}
 
in the middle of the shock layer, we can without loss of generality translate the coordinate $\zeta$ so that $\zeta = 0$ when the travelling wave $\Gamma$ intersects $\Sigma$. It is well-known that an analytic {\it Evans function} can be defined away from the essential spectrum, whose roots (and their multiplicity) coincide with the values in the point spectrum (and their order, respectively); see e.g. \cite{kapprom}. This defines a shooting problem posed with respect to the section $\Sigma$. \\

For sectorial operators it can be shown (see \cite{AGJ,GJ,henry}) that a sufficiently large contour $K \subset \Omega$ can be selected so that the entire point spectrum lies inside $K$. Indeed, we now mimic the argument given in the proofs of Proposition 2.2 and Corollary 2.1 in \cite{AGJ} to show that the operator $\mathcal{L}$ is sectorial. Setting 
$$
Y := -\ve^2 U_{\zeta \zeta} + D(\bar{U})U, 
$$
we can re-write \eqref{eq:linearisedoperator} as 

\begin{equation}\label{eq:bound1}
\begin{aligned}
\ve^2 U_{\zeta \zeta} & = D(\bar{U})U - Y, \\ 
Y_{\zeta \zeta} &= \lambda U - (R'(\bar{U}) + c)U_\zeta.
\end{aligned}
\end{equation}
 
Writing $y := \zeta |\lambda|^\frac14$ and $\tilde{Y}:=Y/|\lambda|^\frac12$, we can rewrite \eqref{eq:bound1} as 

\begin{equation} \label{eq:bound2}
\begin{aligned}
\ve^2 U_{yy}& = \frac{D(\bar{U})U}{|\lambda|^\frac12} - \tilde{Y} \\ 
\tilde{Y}_{yy}& = e^{i \arg \lambda} U - \frac{(R'(\bar{U}) + c)|\lambda|^\frac14 U_y}{|\lambda|}.
\end{aligned}
\end{equation}
As $|\lambda| \to +\infty$, the equation \eqref{eq:bound2} limits to the constant coefficient system

\begin{equation} \label{eq:bound2}
\begin{aligned}
\ve^2 U_{yy}& = - \tilde{Y} \\ 
\tilde{Y}_{yy}& = e^{i \arg \lambda} U.
\end{aligned}
\end{equation}

Writing \eqref{eq:bound2} as a $4\times 4$ system
\begin{equation}
\begin{pmatrix}
U \\  U_y \\ \tilde{Y} \\ \tilde{Y}_y 
\end{pmatrix}_y = 
\begin{pmatrix}
0 & 1/\ve & 0 & 0 \\ 
0& 0& -1/\ve& 0 \\
0 & 0 & 0 & 1 \\ 
e^{i \arg \lambda} & 0 & 0 & 0 
\end{pmatrix}
\begin{pmatrix}
U \\  U_y \\ \tilde{Y} \\ \tilde{Y}_y 
\end{pmatrix}
=: B 
\begin{pmatrix}
U \\  U_y \\ \tilde{Y} \\ \tilde{Y}_y 
\end{pmatrix},
\end{equation}
we find the following (spatial) eigenvalues of $B$:

\begin{equation}\label{eq:boundingeigs}
\mu = \pm \frac{e^{i \frac{\arg \lambda}{4} }}{\sqrt{2 \ve} }\left(1 \pm i\right).
\end{equation}

The key observations here are that if $0<\gamma<\pi/2$ and $|\arg \lambda| < \pi - \gamma$, then there are two (spatial) eigenvalues of $B$ with positive real part and two with negative real part. The corresponding unstable subspaces form an attracting set relative to the flow induced on the Grassmanian $Gr(2,4)$ by equation \eqref{eq:bound2}. Therefore, if $|\lambda| \gg 1$ is sufficiently large, the unstable 2-plane bundle $\varphi^{+}(\zeta,\lambda,\eps)$ stays close to this attractor, and hence can not connect to $\varphi^{-}(\zeta,\lambda,\eps)$ to form a (temporal) eigenvalue. We also observe that the spatial eigenvalues \eqref{eq:boundingeigs} remain uniformly separated as $\ve \to 0$ (in contrast to the example in \cite{AGJ}), and thus we can conclude the following (c.f. Proposition 2.2 and Corollary 2.1 in \cite{AGJ}):

\begin{prop}
There are positive constants $M$ and $\delta$ ($0<\gamma<\pi/2$) so that if $|\lambda|>M$ and $|arg \lambda| < \pi -\gamma$ (with $-\pi <\arg \lambda \leq \pi$), then $\lambda \not \in \sigma(\mathcal{L})$. 
\end{prop} 

We also have the following result asserting the existence of a maximal contour for the point spectrum.
\begin{coro}
There exists a $\beta<0 $ and a simple closed curve $K$, so that $\sigma(\mathcal{L}) \cap \left\{ \lambda: \re \lambda > \beta \right\} $, the spectrum of the operator $\mathcal{L}$ with real part greater than $\beta$, is contained in the interior of $K$. 
\end{coro}

\section{Computing the roots of a Riccati-Evans function} \label{sec:riccatievans}

We now compute the point spectrum of  the operator $\mathcal{L}$. The first step is to formulate a Riccati-Evans function to solve the corresponding shooting problem. Following Ledoux et al. \cite{ledoux2} and Harley et al. \cite{harley}, we think of the flows of $\varphi^{\pm}(\zeta)$ as evolving elements in the Grassmannian $Gr(2,4)$ of complex 2-planes in $\mathbb{C}^4$. The Grassmannian has the structure of a complex manifold (which can be seen via e.g. the well-known {\it Pl\"ucker embedding}). A $2$-plane $P$ in $\mathbb{C}^4$ can be represented by any $4\times 2$ matrix, called a {\it frame} for $P$, whose columns are vectors in $\mathbb{C}^4$ that form a basis of the 2-plane. The Grassmannian has a covering by local coordinate patches $$U_\mathscr{i} = \{Y \in Gr(k,n): Y\cap Y_{i^o} = \emptyset\},$$ where $\mathscr{i} = \{(i_1,\cdots,i_k): 1 \leq i_1 < i_2 < \cdots < i_k \leq n$ is a multi-index and $Y_{i^o}$ denotes the $(n-k)$-plane in $\mathbb{C}^n$ spanned by the vectors $\{e_j\}_{j \notin \mathscr{i}}$. Each element $A \in U_\mathscr{i}$ has a unique matrix representation $A_{\mathscr{i}}$ with respect to the multi-index $\mathscr{i}$ of $A$, with the property that the corresponding $k\times k$ submatrix of $A_{\mathscr{i}}$ is the identity matrix $I_k$. See Sec. 2 of \cite{ledoux2} for details.\\

For concreteness, let us first assume that both subbundles $\varphi^{\pm}(\zeta)$ remain on the coordinate patch $U_{\mathscr{i}}$ with $\mathscr{i} = \{1,\cdots,k\}$ for all $\zeta \in \mathbb{R}$. Then for any $2$-plane $P$ with arbitrarily chosen frame $V$, we arrive at the following (unique) representation on $U_{\mathscr{i}}$ via right multiplication:
\begin{align*}
V = \begin{pmatrix} X(\zeta) \\ Y(\zeta) \end{pmatrix} &\sim \begin{pmatrix} I_k \\ Y(\zeta)X(\zeta)^{-1} \end{pmatrix},
\end{align*}

where we have decomposed the original frame for $P$ into two complex-valued $2\times 2$ matrices $X,\,Y$. Let us now specify the representation of the $2$-plane dynamics on $\mathcal{U}_{\mathscr{i}}$ that is induced by the eigenvalue problem \eqref{eq:linslow}. First perform a block decomposition of the nonautonomous linear operator right-hand side into $2\times 2$ matrices:
\begin{equation}
\label{eq:blockdecomp}
\begin{aligned}
\left(
\begin{array}{cc|cc}
0 & 1/\eps & 0 &0 \\
D(\bar{u})/\eps & 0 & 0 & -1/\eps\\
 \hline 
\lambda - R'(\bar{u}) & 0 & 0 & 0 \\
-c & 0 & 1 & 0
\end{array}
\right)
 &=\left(
    \begin{array}{c|c}
      A & B\\
      \hline
      C & D
    \end{array}
    \right).
\end{aligned}
\end{equation}

Let us write $W(\zeta) = Y(\zeta) X(\zeta)^{-1}$. After some matrix algebra (see \cite{harley}) we arrive at the $2\times 2$ complex-valued {\it matrix Riccati equation}
\begin{equation}
\label{eq:matrixriccati}
\begin{aligned}
W' &= C + DW - WBW.
\end{aligned}
\end{equation} 

Suppose $T$ is an invertible $4\times 4$ matrix representing a change in chart, with $\det(T) = 1$. Letting

\begin{equation}
\label{eq:decompchartchange}
\begin{aligned}
\begin{pmatrix} X_T \\ Y_T \end{pmatrix} &:= T\begin{pmatrix} X\\Y \end{pmatrix} \text{ and}\\
\begin{pmatrix} A_T & B_T \\ C_T & D_T \end{pmatrix} &:= T\begin{pmatrix} A & B \\ C & D \end{pmatrix} T^{-1}
\end{aligned}
\end{equation} 

and defining $W_T = X_T Y_T^{-1}$, we arrive at the corresponding matrix Riccati equation with respect to the new chart:
\begin{equation}
\label{eq:matrixriccati}
\begin{aligned}
W_T' &= C_T + D_T W_T - W_TB_TW_T.
\end{aligned}
\end{equation} 

In this paper we fix the transformation
\begin{align*}
T &= \begin{pmatrix} -i & 0 & 1 & 0 \\ 0 & i & 0 & 1 \\ 0 & 0 & 1 & 0 \\ 0 & 0 & 0 & 1 \end{pmatrix},
\end{align*}

which we selected in an attempt to minimize singularities of the matrix Riccati equation arising from leaving the `natural' chart $U_{\mathscr{i}}$ (note the complex shear transformation $u \mapsto -iu + p$ in the first row). With this choice of chart, we find a simple pole $\lambda_p$ which lies near the real axis. We use this chart for our Riccati-Evans calculations throughout, because it appears to give robust numerical results.\\

Let $z_0 \in \Gamma_{\eps} \cap \Sigma$ denote the (unique) intersection of the travelling wave with the section $\Sigma$ defined in \eqref{eq:sigmasection}. The Riccati-Evans function on this chart is then given by
\begin{align} \label{eq:riccatievans}
E_T(z_0,\lambda)  = \det(W_T^s (z_0;\lambda) - W^u_T(z_0;\lambda)),
\end{align}

where $W_T^{u,s}$ denote the corresponding values of $\varphi^{+,-}$, respectively at the point of intersection (and written in the same chart coordinates). Using the fact that $\det(T) = 1$ and the argument principle, we can assess the roots and poles of the Riccati-Evans function (and hence deduce the zeroes of the Evans function) by computing its winding number along judiciously chosen contours to the right of the essential spectrum.

\begin{figure}[t!] 
\centering
(a)\includegraphics[width=0.85\textwidth]{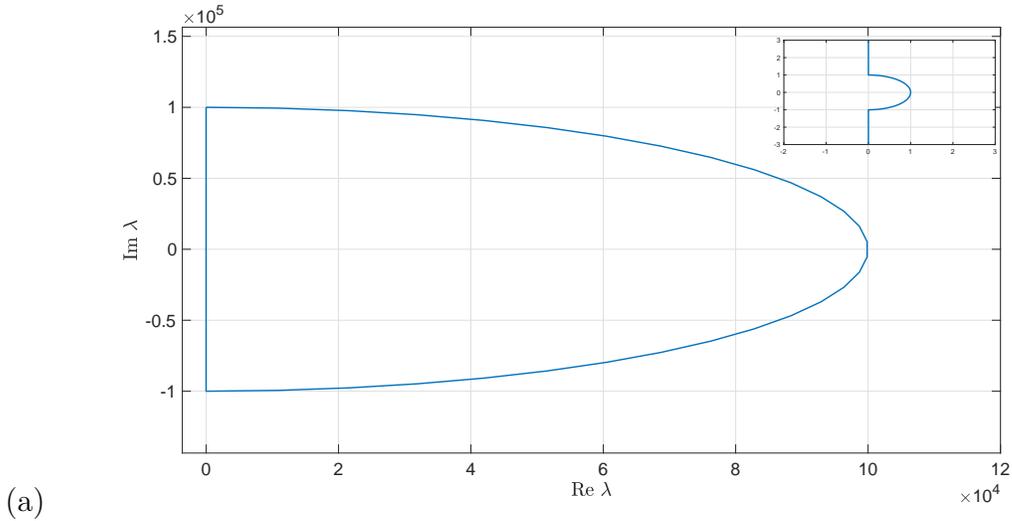}\\
(b)\includegraphics[width=0.85\textwidth]{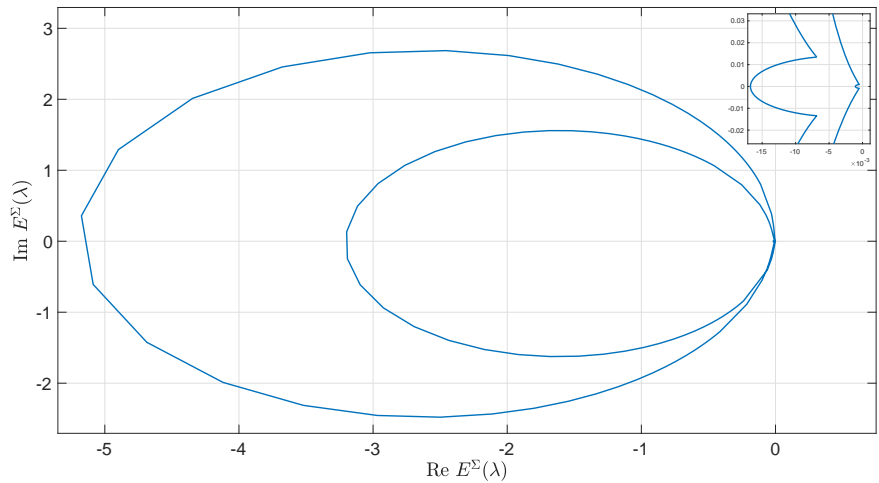}
\caption{Computations of the winding of the Riccati-Evans function \eqref{eq:riccatievans} for the parameter set $(\eps,c) = (10^{-4},0.19686)$.  (a) Semicircular contour in the right-half complex plane of radius $10^5$ centered at the origin, with its diameter aligning with the imaginary axis and also containing a small semicircular detour of radius 1 to avoid the translational eigenvalue at $\lambda = 0$ (see inset). (b) Evaluation of the Riccati-Evans function along this contour. The image is a closed curve that does not wind around the origin (see inset)---thus, the Riccati-Evans function has a winding number of 0 along the chosen contour.}
\label{fig:riccatievans}
\end{figure}

\begin{figure}[t!] 
\centering
(a)\includegraphics[width=0.85\textwidth]{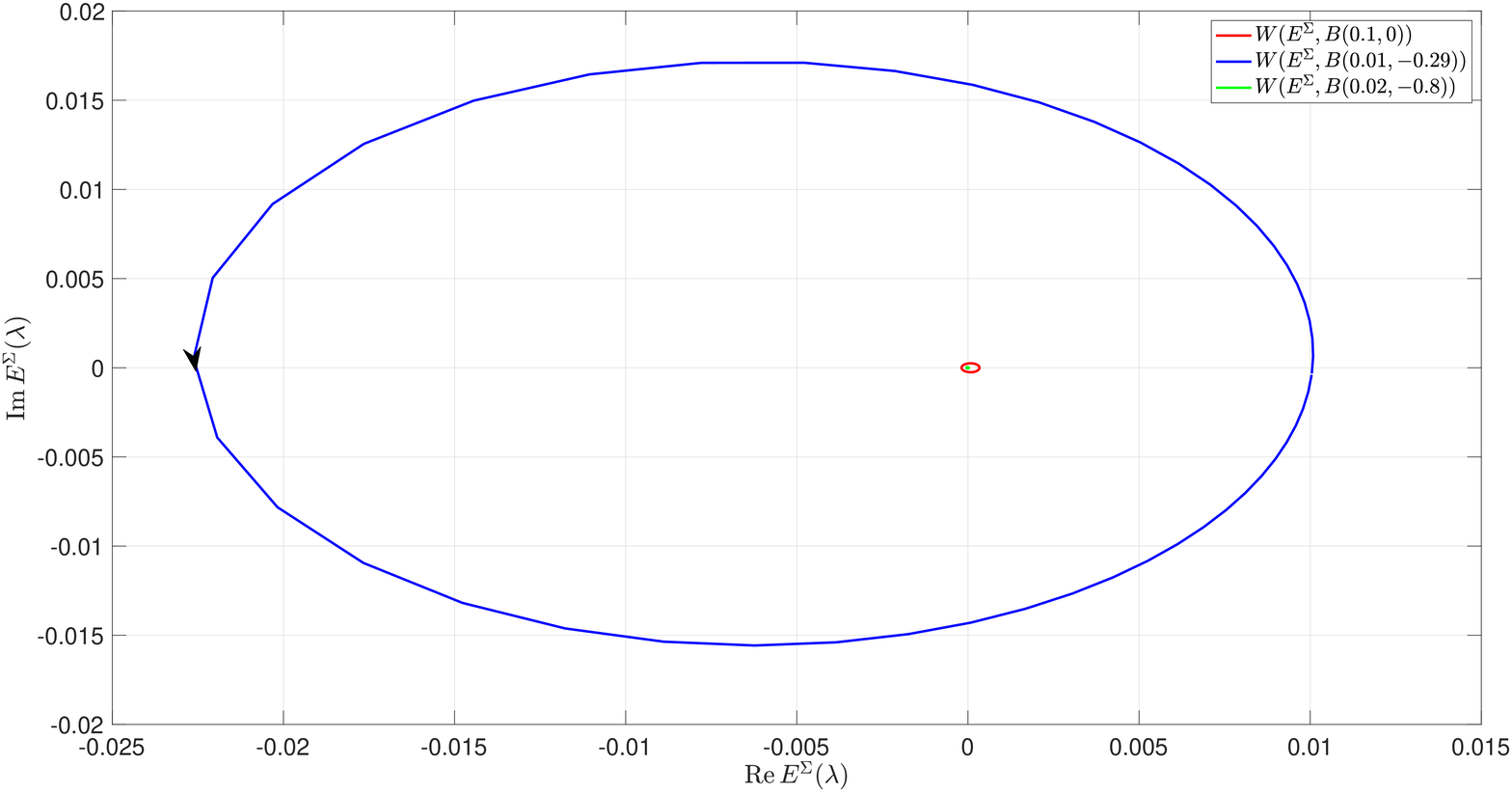}\\
(b)\includegraphics[width=0.85\textwidth]{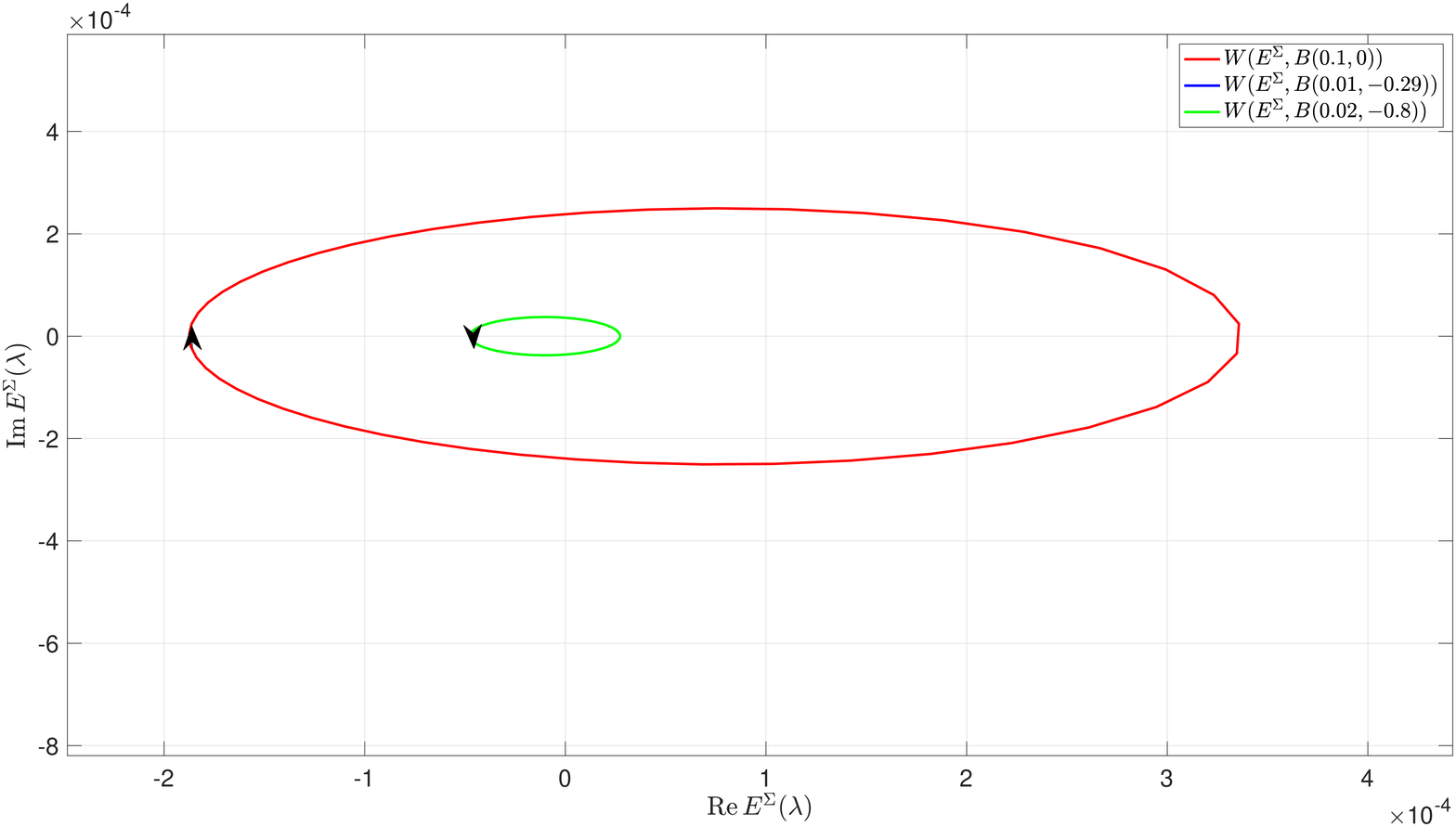}
\caption{Computations of the winding of the Riccati-Evans function along the circle (a) $B(0.1,0)$ as well as along  (b) $B(0.01,-0.29)$ and $B(0.02,-0.8)$. The corresponding windings are $1$, $-1$, and $1$, respectively, which is consistent with the existence of roots at $\lambda_0 = 0$ and $\lambda_1 \approx -0.29$ and a pole at $\lambda_p = -0.29$.}
\label{fig:riccatievans2}
\end{figure}

\begin{remk}
Numerical integration of the matrix Riccati equation \eqref{eq:matrixriccati} appears to be very stable relative to working directly with the linearised system \eqref{eq:linslow} or its induced two-plane dynamics via derivations (see e.g. \cite{jonestin}). The simpler case of tracking an unstable {\it line} bundle is a useful starting point to understand this numerical robustness. The obvious strategy there is to consider the projectivization of the eigenvalue problem and to study the dynamics in projective space over a suitable chart. Then the `most unstable direction' becomes a curve of {\it attracting fixed points} over the slow manifolds (with respect to the corresponding family of `frozen' systems associated with the projectivization of \eqref{eq:linslow}, where the time variable $\zeta$ in the nonautonomous component is frozen).  See e.g. \cite{GJ,lizarraga}, where this observation is crucial in the construction of so-called {\it relatively invariant attracting sets} over the strong unstable direction. As a heuristic, attracting sets tend to be easier to numerically approximate.\\

Turning to the present case, we observe numerically that the unstable $2$-plane bundle similarly becomes an attracting set with respect to the matrix Riccati equation \eqref{eq:matrixriccati}.  Let us recall that the induced derivation $$w' = A^{(k)}w,$$ specifying the $k$-plane dynamics of a linear system $z' = Az$, has eigenvectors corresponding to $k$-fold wedge products of the eigenvectors of the original system, with the corresponding eigenvalues given by the $k$-fold sums of the appropriate eigenvalues of the original system (see e.g. \cite{bourbaki})---in particular, the unstable subspace formed by taking the wedge product of two unstable eigenvectors corresponds to an eigenvector of the new induced system, corresponding to the new largest positive eigenvalue.\\

  Then our present numerical observations are consistent with the fact that at the projective level, these `unstable eigenplanes' become attracting fixed points, in analogy with the line bundle case. The numerical integration of \eqref{eq:matrixriccati} is also stable when integrating the stable subbundle backwards in time, for the same reason. We expect that the matrix-Riccati equation can be a useful starting point to write down the corresponding estimates for relatively invariant attracting sets over unstable $k$-plane bundles in $\mathbb{C}^n$. \end{remk}

As shown in Fig. \ref{fig:riccatievans}, an evaluation of the winding number of the Riccati-Evans function \eqref{eq:riccatievans} on a large contour in the right-half plane shows that there appear to be no eigenvalues of positive real part. Our function appears to map increasingly large semicircular contours onto a compact set bounded in the left-half complex plane, and thus it does not appear possible to generate a winding around the origin by choosing larger contours with the same shape. \\

We also compute a counterclockwise winding of $+1$ along large circular contours containing the origin and to the right of the essential spectrum. In view of the calculation in Fig. \ref{fig:riccatievans}, we decide to investigate the winding along much smaller circles centered on the real line near the origin. As shown in Fig. \ref{fig:riccatievans2}, we evaluate the winding of the Riccati-Evans function about two small circles $\lambda_0 = 0$ and $\lambda_1 \approx -0.8$ to find an index of $+1$, which is consistent with the existence of simple roots near these values. We also find a pole $\lambda_p$  of order 1 (i.e. having a winding of $-1$), with $\text{Re}(\lambda_1) < \text{Re}(\lambda_p) < \text{Re}(\lambda_0)$. The sum of the winding numbers of $\lambda_0$, $\lambda_p$, and $\lambda_1$ gives a total index of $+1$, which is consistent with the winding number evaluated about the large circles.\\

Further investigations along different choices of contour did not reveal the existence of any other roots or poles. We conclude that the point spectrum consists only of the simple eigenvalues $\lambda_0,\,\lambda_1$ with
\begin{equation}
\label{eq:fulleigs}
\begin{aligned}
\lambda_0 &= 0 \text{ and}\\
\text{Re}(\lambda_1) &\approx -0.8.
\end{aligned}
\end{equation}

Thus, the corresponding shock-fronted travelling wave is nonlinearly stable for the singular perturbation value $\eps = 10^{-4}$. Investigations of a few different small values of $\eps$ ranging from $10^{-4}$ to $10^{-2}$  gave very similar numerical results.

\section{Fast (un)stable-to-(un)stable connections} \label{sec:fastconnections}

We now show how the results in the previous section can be explained via a fast-slow splitting of the (un)stable $2$-plane bundles into fast and slow line subbundles. Such splittings were analyzed in great depth in the papers \cite{AGJ,GJ}; we relate the geometric construction here. The fundamental observation is the topological characterization of the point spectrum within a contour via the first Chern class of a complex vector bundle over a sphere called the {\it augmented unstable bundle}. The sphere is formed by extending the contour along the time direction to define an infinite cylinder, and then applying a compactification; this sphere is endowed with the $2$-plane unstable bundle defined using the flow of the eigenvalue problem.\\

 The first Chern number is robust to small deformations of the complex vector bundle and also distributes over Whitney sums of line bundles \cite{atiyah}. The latter property potentially allows one to compute the windings of a $2$-plane bundle along the contour $K$ by summing the windings of fast and slow line bundles along $K$, and the former potentially allows a straightforward computation of these `sub-windings' for sufficiently small values of $\eps>0$, by means of {\it reduced} eigenvalue problems. \\
 
 Our first task is demonstrate the existence of {\it fast unstable-to-unstable} connections of line bundles for open sets of $\lambda \in \Omega$ and for $\eps > 0$ sufficiently small. Recall the hyperbolic fast-slow eigenvalue splittings $\mu_i^{\pm}$ of the eigenvalue problem near the tails of the wave. We say that there is a fast unstable-to-unstable connection for $\lambda$ if the {\it fast unstable line bundle} $\hat{e}_1$, defined as the unique  bundle that converges to the eigenspace of the eigenvalue of largest real part $\text{Re}(\mu^+_1)>0$ of $a_+(\lambda,\eps)$ as $\zeta \to -\infty$ under the flow of the eigenvalue problem, converges to the eigenspace of the eigenvalue of largest real part $\text{Re}(\mu^-_1) > 0$ of $a_-(\lambda,\eps)$ as $\zeta \to +\infty$.  We can analogously define a {\it fast stable-to-stable connection} $\hat{e}_4$ by reversing time and asking for connections from the eigenspace of the smallest eigenvalue $\mu^-_4 < 0$ for $A_-(\lambda,\eps)$ to that corresponding to $\mu^+_4<0$ for $A_+(\lambda,\eps)$, respectively.\\
 
 An important reduced system in this context is the {\it fast reduced eigenvalue problem} defined from a singular limit of \eqref{eq:linfast}, i.e.
 
 \begin{equation}
\label{eq:reducedfasteig}
\begin{aligned}
u' &= w\\
w' &= D(\bar{u})u.
\end{aligned}
\end{equation} 

This singular limit is defined along the (singular) shock-front of the wave.  When the fast reduced eigenvalue problem has a fast unstable-to-unstable connection for $\lambda \in \Omega$, a key estimate known as the {\it elephant trunk lemma} (see \cite{GJ}) allows this connection to be extended along the entire wave for small values of $\eps>0 $. Unfortunately, since $\lambda$ does not appear in our fast reduced eigenvalue problem, it degenerates so that there is always a fast unstable-to-{\it stable} connection in the singular limit!\footnote{This is true because \eqref{eq:reducedfasteig} is the variational equation of the layer vector field \eqref{eq:layerproblem} along the singular shock front connecting two saddle points, and therefore automatically has a nontrivial solution that vanishes exponentially quickly as $\xi \to \pm \infty$.} This motivates a numerical analysis to test whether this connection indeed breaks into the desired fast unstable-to-unstable connection when $\eps > 0$.\\

\begin{figure}[t!] 
\centering
(a) \includegraphics[width=0.95\textwidth]{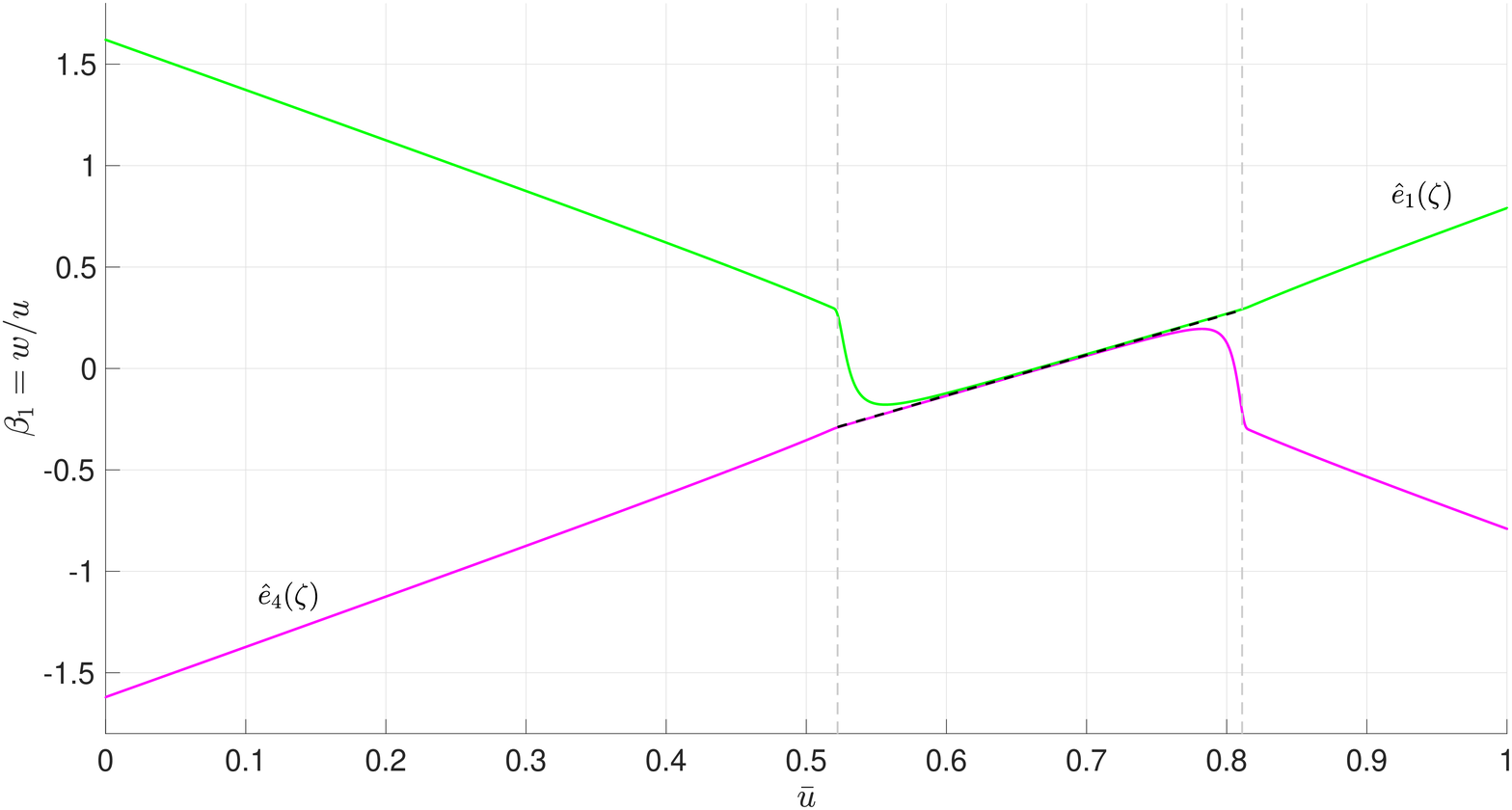}
\caption{A depiction of the fast unstable-to-unstable connection (green) and a fast stable-to-stable connection (magenta) when $(\eps,c,\lambda) = (10^{-4},0.19686,0)$, using the projectivization $\beta_1 = w/u$. The gray dashed lines denote the jump values $\bar{u}=\bar{u}_{\pm}$ and the black dashed line segment depicts the singular fast unstable-to-stable connection. Different choices of $\lambda$ give very similar graphs.}
\label{fig:nearmiss}
\end{figure}

In Fig. \ref{fig:nearmiss} we depict an example of such fast (un)stable-to-(un)stable connections when $\lambda = 0$. We choose a chart of projective space $\mathbb{CP}^3$ with $(\beta_1,\beta_2,\beta_3) = (w/u,p/u,v/u)$ to depict the dynamics, because in the singular limit of \eqref{eq:linslow} we can write down the exact expressions 
\begin{align*}
v^s(\bar{u}) = \begin{pmatrix} \sqrt{D(\bar{u})} \\ -D(\bar{u})\\0\\0 \end{pmatrix} &\text{ and } v^u(\bar{u}) = \begin{pmatrix} \sqrt{D(\bar{u})} \\ D(\bar{u})\\0\\0\end{pmatrix}
\end{align*}

characterizing the fast stable and unstable eigendirections for the eigenvalues $\mu_4(\bar{u}) = -\sqrt{D(\bar{u})}$ and $\mu_1(\bar{u}) = +\sqrt{D(\bar{u})}$, respectively. With respect to the chart representation, we have $v^s(\bar{u}) = (-\sqrt{D(\bar{u})},0,0)$ and  $v^u(\bar{u}) = (\sqrt{D(\bar{u})},0,0)$. This accounts for the nearly linear dynamics away from the jump in Fig. \ref{fig:nearmiss}. We can also read off the kind of connection made by checking whether the trajectory stays above $\beta_1 > 0$ (unstable-to-unstable) or below $\beta_1 < 0$ (stable-to-stable) away from the shock.

\begin{remk}
Recall that the case $\lambda = 0$ corresponds to the variational problem, where we already know there is an eigenvalue! The relevant connection in this case occurs within the slow subbundle instead (see the discussion below \eqref{eq:linslow}). We discuss the relevant slow eigenvalue problem in the next section.
\end{remk}

 We test the possibility of a fast unstable-to-stable connection (i.e. a {\it fast eigenvalue}) by defining a proxy {\it fast reduced Riccati-Evans function} for $\eps > 0$ on the section $\Sigma$ as follows. Because we consider line bundles in $\mathbb{C}^3$, our strategy is to consider a projectivisation of the system \eqref{eq:linslow} on a chart with $u \neq 0$, with $(\beta_1,\beta_2,\beta_3) = (w/u,p/u,v/u)$, and to define
 \begin{equation}
\label{eq:fastevans}
\begin{aligned}
E_f^{\Sigma}(\lambda) &:= \beta_1^{u}(z_0,\lambda) - \beta_1^{s}(z_0,\lambda),
\end{aligned}
\end{equation} 
where $z_0$ is defined on the section $\Sigma$ as before. This function is `Evans-like' in the sense that if there are no roots of this function within a contour, then it is hopeless for a fast unstable-to-stable connection to be made for values of $\lambda$ within the contour: such a connection is necessarily unique. We highlight, however, that a zero of $E_f^{\Sigma}$ does {\it not} necessarily imply the existence of a fast eigenvalue: the values of $\beta_2,\,\beta_3$ must also coincide.  Regardless, it is meaningful to test the winding number of $E_f^{\Sigma}(\lambda)$ evaluated along contours in $\mathbb{C}$ in order to efficiently test candidates $\lambda$ for fast connections.\\

We found only one simple zero of $E_f^{\Sigma}$ at the value $\lambda \approx 3718.025$. However, at this value it can be computed directly that $\beta_2^{u}(z_0,\lambda) - \beta_2^{s}(z_0,\lambda) \approx -5.08$, and thus there is still a fast unstable-to-unstable connection made at this value. Using \eqref{eq:riccatievans}, we can verify that there is indeed no point spectrum in this parameter region.\\

We found that the additional fast direction endowed by nonlocal regularization (compared to viscous relaxation) appears to give us only trivial fast unstable (resp. fast stable) connections, and thus the winding of the Evans function for this two-plane bundle with this extra fast direction is not the same as that without it, and this extra fast direction does not contribute to the point spectrum.

\section{The slow eigenvalue problem} \label{sec:sloweigs}

Having shown that the fast connections do not contribute to the point spectrum, we now demonstrate that a {\it slow eigenvalue problem}, defined near the tails of the shock-fronted travelling wave, controls the point spectrum of the full eigenvalue problem when $\eps > 0$ is sufficiently small. In analogy to the definition of \eqref{eq:reducedfasteig}, we take the singular limit $\eps \to 0$ of \eqref{eq:linfast} to obtain the {\it reduced slow eigenvalue problem}
 \begin{equation}
\label{eq:reducedsloweig}
\begin{aligned}
\dot{P}&= (\lambda-R'(\bar{U}))(V/D(\bar{U}))\\
\dot{V} &= P - c_0 V/D(\bar{U})
\end{aligned}
\end{equation} 

defined with respect to the linearised constraint $V = D(\bar{U})U$. The reduced slow eigenvalue problem is defined on the outer saddle-type branches of the critical manifold $S_{0}$, in particular along the singular limit of the travelling wave lying on these branches (i.e., the left and right slow segments of the red trajectory depicted in Fig. \ref{fig:transversehet}). Reduced eigenvalue problems of this sort appear already in \cite{GJ}. We also point the reader to the derivation in the viscous case \cite{lizarraga}: we note that the reduced problem in the present nonlocal case is topologically equivalent (modulo desingularisation and a trivial orientation reversal in the `time' variable $\zeta$).\\

For $\lambda \in \Omega$, the linear operator for the reduced slow eigenvalue problem remains hyperbolic near the end states $\bar{U} = 0$ and $\bar{U} = 1$, and so we retain a well-defined geometrical characterization of the eigenvalue problem in terms of finding an unstable-to-stable connection. Specifically, we can track the reduced slow unstable line bundle $\varphi^{+}_s$ from $\bar{U} = 1$ and check whether it connects with the reduced slow stable  line bundle $\varphi^{-}_s$ emanating from $\bar{U} = 0$. \\

There is one technical issue to settle: since the reduced slow eigenvalue problem is defined on disjoint segments of the critical manifold, we must also define a {\it jump map} which specifies how data at the beginning of the shock on $u=u_+$ is transported to $u=u_-$ at the end of the shock. This map is derived as a singular limit in Sec. 7.3 of  \cite{lizarraga}: via an identical procedure, we obtain the jump map $J_{\lambda}:\mathbb{C}^2 \to \mathbb{C}^2$ given by
\begin{align} \label{eq:jumpmap}
J_{\lambda}\begin{pmatrix} P \\ V \end{pmatrix} &= \begin{pmatrix}1 & \frac{R(\bar{U})_+-R(\bar{U}_-) - \lambda(\bar{U}_+ - \bar{U}_-)}{cU_{+} - \bar{P}_+} \\ 0 & \frac{c \bar{U}_+ - \bar{P}_+}{cU_{-}-\bar{P}_-}  \end{pmatrix} \begin{pmatrix} P \\ V \end{pmatrix}.
\end{align}

\begin{figure}[t!] 
\centering
(a)\includegraphics[width=0.8\textwidth]{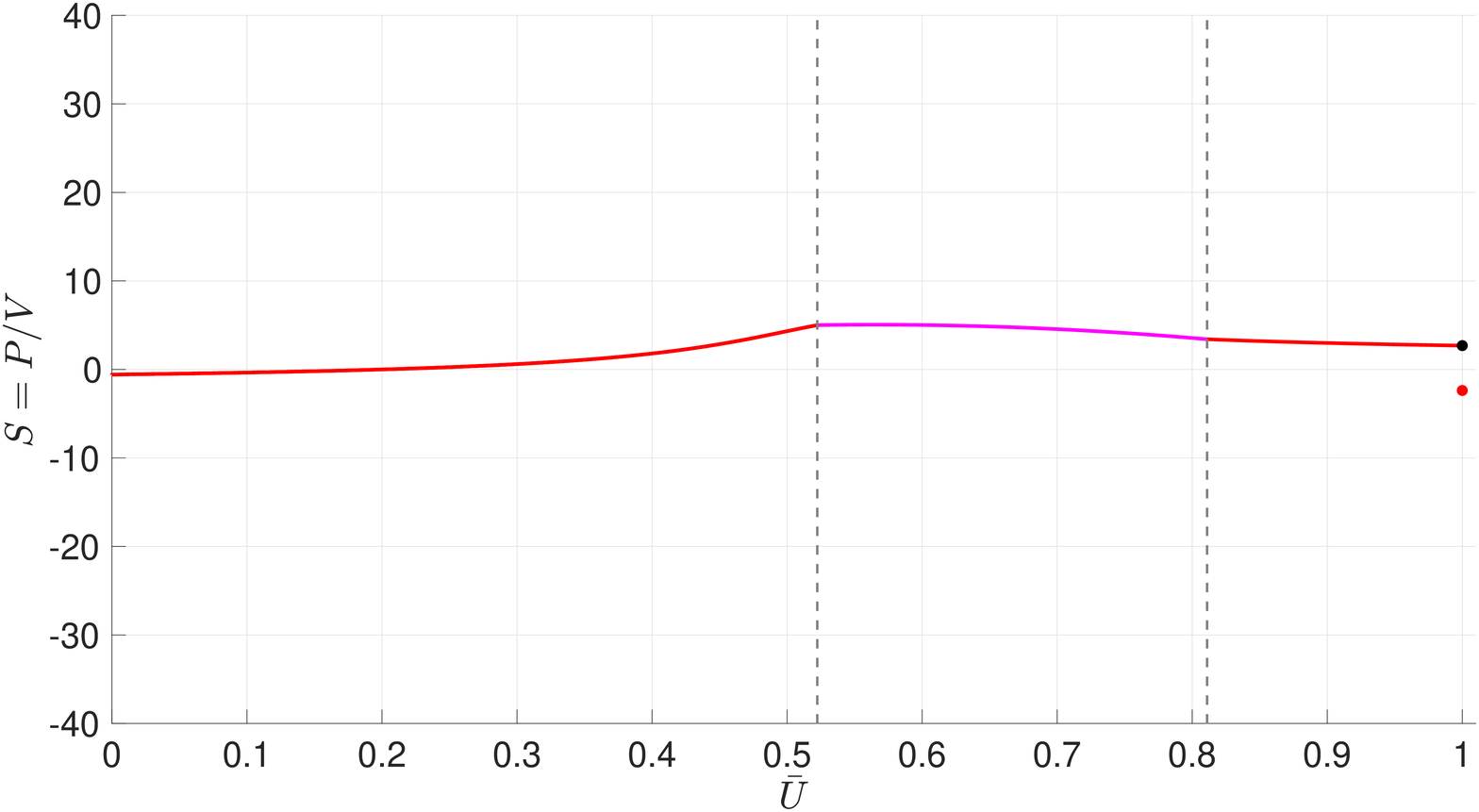}
(b)\includegraphics[width=0.8\textwidth]{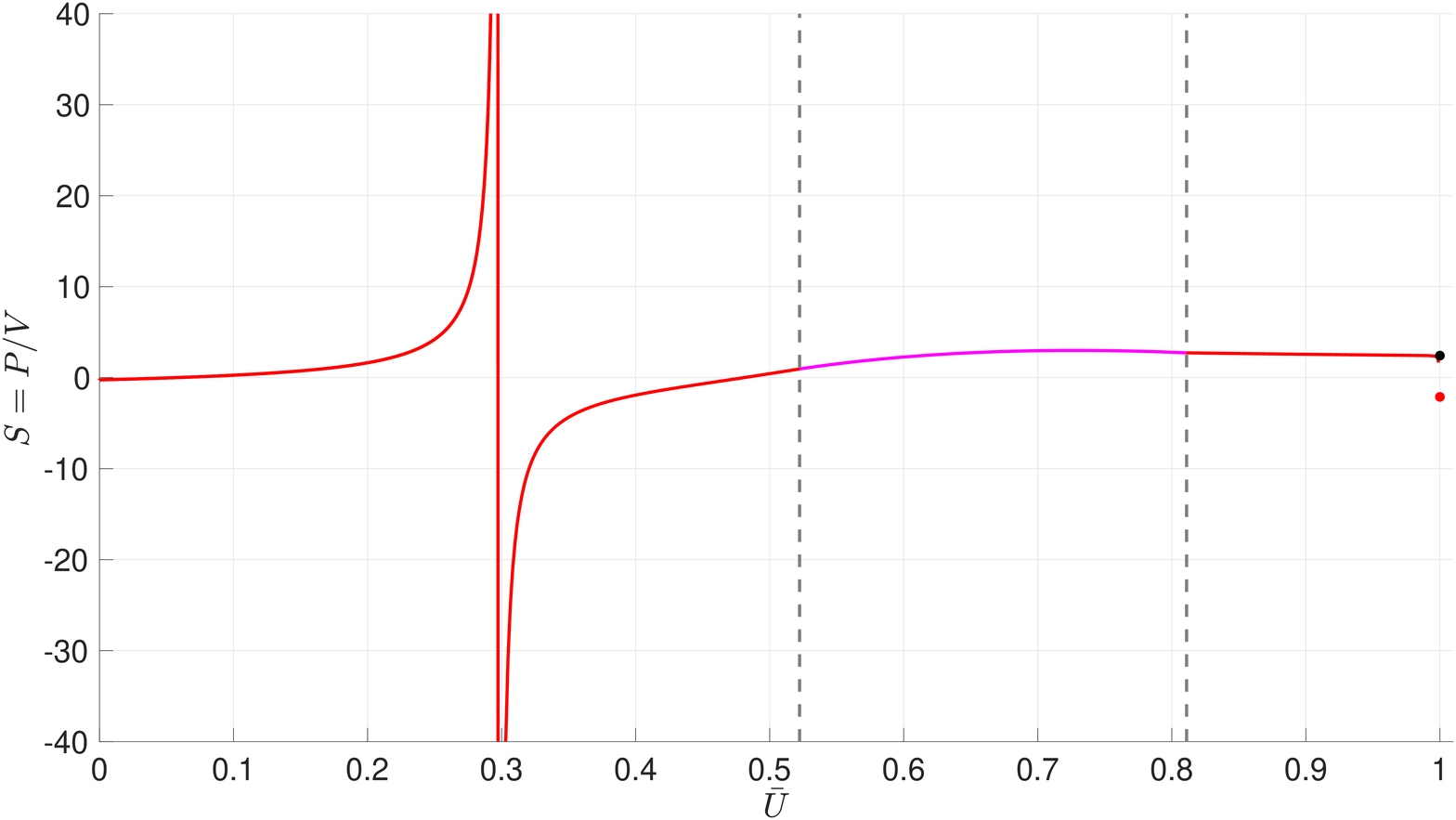}
\caption{Slow eigenvalues (fast unstable-to-stable connections) of the reduced slow eigenvalue problem \eqref{eq:reducedsloweig}--\eqref{eq:jumpmap} when (a) $\lambda_0 = 0$ and (b) $\lambda_1 \approx -0.80031$. Red point: attracting fixed point corresponding to projectivized unstable direction at $\bar{U} = 1$; black point: saddle fixed point corresponding to projectivized stable direction at $\bar{U} = 1$. The red segments depict the flow of \eqref{eq:reducedsloweig} and the intermediate magenta segments depict the output of jump map \eqref{eq:jumpmap}.}
\label{fig:sloweigs}
\end{figure}

Note that the jump map is linear, which implies that we can write down a well-defined projectivization of the map. In the framework of \cite{AGJ,GJ}, the jump map serves as a {\it clutching function} over the equator of the sphere when defining the reduced slow line bundle. We also highlight that the jump map is defined for $\lambda = 0$, i.e. it is compatible with the evaluation of the reduced vector field \eqref{eq:reduced2} at either end of the jump point.\\

\begin{figure}[t!] 
\centering
\includegraphics[width=0.95\textwidth]{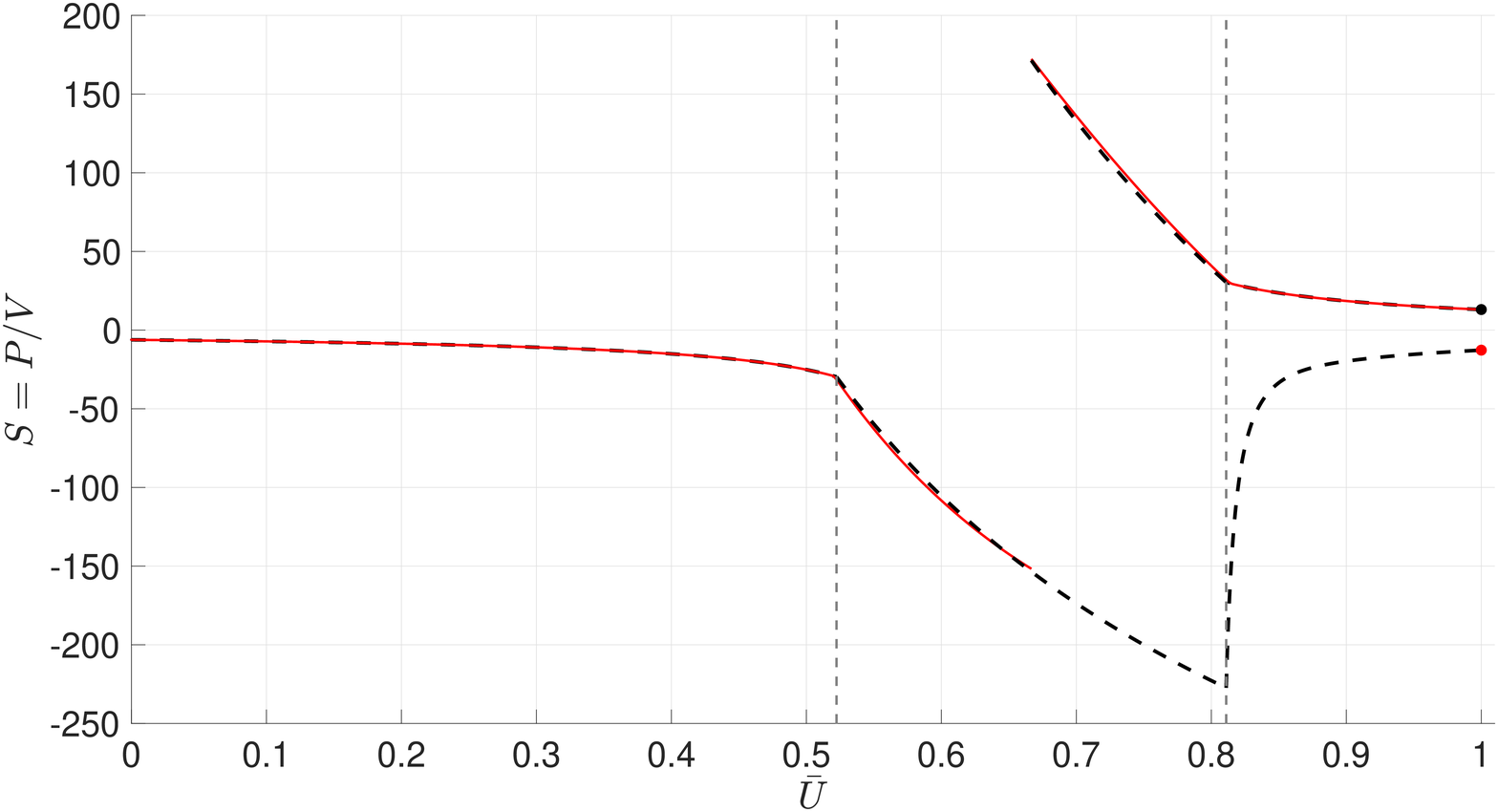}
\caption{Demonstration of a {\it slow} unstable-to-unstable connection for $\lambda = 100$. Red curve segments: solutions of the boundary-value problems beginning on the weak unstable eigendirection (i.e. corresponding to the eigenspace of $\mu^{+}_2$) for the saddle-point at $\bar{U} = 1$ with the section $\Sigma = \{\bar{U} = 2/3\}$, and beginning on the weak stable eigendirection (i.e. corresponding to the eigenspace of $\mu^{-}_3$ for the saddle-point at $\bar{U} = 0$) with that section. Red point: attracting fixed point corresponding to projectivized unstable direction at $\bar{U} = 1$; black point: saddle fixed point corresponding to projectivized stable direction at $\bar{U} = 1$. Black dashed segments: solutions of the projectivization of \eqref{eq:reducedsloweig} on the chart $S = P/V$ (with $V \neq 0$) together with output from the jump map \eqref{eq:jumpmap}.  The red segments were computed using the \texttt{bvp4c} boundary-value solver package in MATLAB, with a relative error tolerance of $10^{-8}$.}
\label{fig:slowconnection}
\end{figure}

Via an identical analysis to that in Sec. 8 of \cite{lizarraga}, we may conclude that the reduced slow eigenvalue problem \eqref{eq:reducedsloweig} together with the jump map \eqref{eq:jumpmap}: (1) does not have eigenvalues with nonzero imaginary part, and (2) does not have any eigenvalues of positive real part. Because of these facts, we may restrict our attention to the one-dimensional problem defined by projectivization of \eqref{eq:reducedsloweig} (by e.g. taking the chart $S = P/V$ with $V\neq 0$) and restricting the parameter set to only real $\lambda$.\\

 By direct computation (see Fig. \ref{fig:sloweigs}) we find only two eigenvalues: 
 \begin{equation}
\label{eq:sloweigs}
\begin{aligned}
\lambda_0 &= 0 \text{ and}\\
\lambda_1 &\approx -0.80031.
\end{aligned}
\end{equation} 
By the analysis in Sec. 8 of \cite{lizarraga}, we may immediately conclude that both eigenvalues are simple. Alternatively, we can also define a section, say $\Sigma_R = \{\bar{U} = 0.7\}$, and define a {\it slow reduced Evans function}
\begin{align*}
E^{\Sigma_R}_s(\lambda) &:= S^+(z_0,\lambda)-S^-(z_0,\lambda),
\end{align*}
 
 where $S^{+}(z_0,\lambda)$  denotes the location of the unstable slow subbundle $\varphi^{+}$ at the (first) intersection $z_0$ of the singular limit of the travelling wave with $\Sigma_R$ (and analogously for $S^-(z_0,\lambda)$). These computations were carried out in Sec. 8 of \cite{lizarraga} for slightly different jump conditions and wavespeed, so we exclude them here while pointing out that the numerical output is equivalent. \\
 
We emphasize that the slow eigenvalue count \eqref{eq:sloweigs} is in excellent agreement with the eigenvalue count of the `full' problem \eqref{eq:fulleigs} which was obtained with Riccati-Evans computations. In fact, we can show even more: by using a boundary value problem solver, we can set up a matching problem for the {\it slow} line subbundles defined by the system \eqref{eq:linslow} on the section $\Sigma$ and show that they are in excellent agreement with the reduced dynamics \eqref{eq:reducedsloweig}--\eqref{eq:jumpmap}. We depict one such computation in Fig. \ref{fig:slowconnection}, where we choose a real value $\lambda = 100$ for convenience.  \\

We conclude our analysis by outlining the argument which relates the point spectrum of the full problem with the point spectra of the fast and slow reduced problems, according to the topological characterization given in \cite{AGJ,GJ}. Fixing a large enough contour $K$ so that the entire point spectrum lies inside it, the associated {\it augmented unstable bundle} of the full problem, which we denote $\mathcal{E}_{\eps}(K)$, can be split into a Whitney sum of fast and slow line subbundles as
\begin{align}
\mathcal{E}_{\eps}(K) &= \mathcal{E}^f_{\eps}(K) \oplus \mathcal{E}^s_{\eps}(K)
\end{align}

when $\eps > 0$ is sufficiently small. It is well-known (see \cite{atiyah}) that the first Chern number evaluation distributes over line bundle decompositions as
\begin{align*}
c_1(\mathcal{E}_{\eps}(K)) &= c_1(\mathcal{E}^f_{\eps}(K)) + c_1(\mathcal{E}^s_{\eps}(K)),
\end{align*}

and by a theorem of Alexander, Gardner, and Jones \cite{AGJ}, the first Chern number of the augmented unstable bundle is equal to the eigenvalue count (with multiplicity) of the linearised operator inside $K$. We can now focus on computing separately the first Chern number of the fast and slow line bundles.\\

Our numerical observations (e.g. Fig. \ref{fig:slowconnection})  suggest the existence of a deformation
\begin{align*}
\mathcal{E}^s_{\eps}(K) &\cong \mathcal{E}^s_{0}(K),
\end{align*}

for $\eps > 0$ sufficiently small, where $\mathcal{E}^s_{0}(K)$ is a reduced augmented line bundle defined by using the reduced eigenvalue problem \eqref{eq:reducedsloweig} together with a clutching operation defined from the jump map \eqref{eq:jumpmap}; see rigorous examples of such constructions in \cite{GJ,lizarraga}. The resulting {\it slow eigenvalues} can be counted as described in this section, i.e. $c_1(\mathcal{E}^s_{\eps}(D)) \leq 2$ depending on the contour $D$ chosen. \\

It remains to evaluate $c_1(\mathcal{E}^f_{\eps}(K))$ for sufficiently small values of $\eps > 0$, and here we observe a surprising departure from the techniques developed in \cite{AGJ,GJ,jones}. We can ask whether there exists a similar deformation onto some singular fast line bundle $\mathcal{E}^f_0(K)$, i.e.
\begin{align*}
\mathcal{E}^f_{\eps}(K) &\cong \mathcal{E}^f_{0}(K).
\end{align*}
In \cite{GJ}, the singular fast line bundle is defined by finding singular fast unstable-to-unstable connections in the fast reduced eigenvalue problem, and then an elephant trunk lemma is proven to show that these connections persist along the entire wave when $\eps > 0$. But our fast reduced eigenvalue problem \eqref{eq:reducedfasteig} degenerates! Nonetheless, it appears to be the case that $c_1(\mathcal{E}^f_{\eps}(K)) = 0$ for all sufficiently small $\eps > 0$, in view of our analysis in Sec. \ref{sec:fastconnections}. \\

This completes our analysis of the fourth-order nonlocal problem: we have $$(c_1(\mathcal{E}_{\eps}(K)) =2$$

for all $\eps > 0$ sufficiently small, where the first Chern number is exactly accounted for by the simple translational eigenvalue at the origin, together with a nontrivial eigenvalue $\lambda_1$ with $-1 < \text{Re}(\lambda_1) < 0$. This nontrivial eigenvalue can be furthermore be well-approximated by the corresponding slow eigenvalue $\lambda_1$ shown in Fig. \ref{fig:sloweigs}. Thus, the family of travelling waves $\{\Gamma(c(\eps),\eps)\}_{\eps \in (0,\bar{\eps}]}$ is nonlinearly stable for each sufficiently small $\eps > 0$.

\section{Concluding remarks} \label{sec:conclusion}

We have given a comprehensive analysis of the nonlinear (asymptotic) stability of a family of shock-fronted travelling waves arising in the regularized RND PDE \eqref{eq:master}. This result is shown by verifying that the linearised operator is sectorial Sec. \ref{sec:eigenvalueproblem}, and then calculating the point spectrum using a {\it Riccati-Evans function} (Sec. \ref{sec:riccatievans}). A fast-slow splitting of the corresponding eigenvalue problem was also furnished (Secs. \ref{sec:fastconnections}--\ref{sec:sloweigs}), allowing us to explain the emergence of each eigenvalue and its order.\\

It would be interesting to analyze such fast-slow splittings in more general models that incorporate both nonlocal {\it and} viscous regularizations as described in \cite{proceedings22}, and in particular whether to assess whether the nonlinear stability of the wave persists for all combinations of these two regularizations. The corresponding stability result for the viscous relaxation limit in \cite{lizarraga} is slightly weaker because of the asymptotically vertical nature of the essential spectrum; on the other hand, there are parameter regimes in the generalised model for which viscous shocks can potentially form, and yet the linearised operator is still sectorial. It is also of interest to determine to what extent slow eigenvalue problems are responsible for the generation of eigenvalues in the point spectrum for other more exotic choices of regularization. \\

In view of the eigenvalue computation using first Chern numbers in the previous section, further analysis (e.g. in the style of \cite{GJ,lizarraga}) that explains the behavior of $\mathcal{E}_{\eps}^f(K)$ for small values of $\eps$ is desirable. We suggest that a trivial `singular fast bundle' $\mathcal{E}^f_0(K)$ should exist, but it should be constructed on the $\zeta$-timescale (and not on the $\xi$-timescale as is done in \cite{GJ}), i.e. it should be constructed using a clutching procedure, similarly to what is done for the singular slow bundle. A rigorous approximation theorem for the fast bundle (i.e. a proof that $\mathcal{E}^f_{\eps}(K) \cong \mathcal{E}^f_0(K)$ for $\eps > 0$ sufficiently small) will require new estimates that do not use elephant trunk lemma-style theory. We suggest that new estimates of exchange lemma-type (see \cite{lizarraga}) may prove useful; the standard $(k+\sigma)$--exchange lemma is ideally suited to tracking plane bundles near the center-unstable manifolds of saddle-type slow manifolds (see e.g. \cite{joneskopell,jonestin}). Melnikov-type analysis might also be useful to describe transversal breaking of degenerate unstable-to-stable connections; see Lemma 1 and Theorem 1, part (ii) in \cite{szmolyan} for such a result in the case of planar viscous shockwaves.

\section*{Acknowledgements}

Both authors acknowledge support by the Australian Research Council Discovery Project grant DP200102130. We would like to thank Bj\"{o}rn de Rijk for a useful discussion about analytical factorizations of Evans functions.

\end{document}